\documentclass[11pt]{article}

\input{amssym.def}
\input{amssym}
\usepackage{eufrak}

\newtheorem{theorem}{Theorem}

\newtheorem{corollary}[theorem]{Corollary}
\newtheorem{definition}[theorem]{Definition}

\newtheorem{lemma}[theorem]{Lemma}
\newtheorem{proposition}[theorem]{Proposition}
\newtheorem{remark}[theorem]{Remark}

\def\Uq{U_q(sl(m))}
\def\gg{\mbox{$\frak g$}}
\def\ot{\otimes}
\def\h{{\hbar}}
\def\K{{\Bbb K}}
\def\C{{\Bbb C}}
\def\R{{\Bbb R}}
\def\Sym{{\rm Sym\, }}
\def\ad{{\rm ad\, }}
\def\End{{\rm End}}
\def\Ren{R_{\End}}
\def\vv{V^{\otimes 2}}

\def\lhq{\ifmmode {\cal L}(q,\hbar)\else ${\cal L}(q,\hbar)$\fi}
\def\lrqh{{\cal L}(R_q,\h)}
\def\lrq{{\cal L}(R_q)}
\def\lrqo{{\cal L}(R_q,1)}
\def\lqh{\ifmmode {\cal L}(q,\hbar)\else ${\cal L}(q,\hbar)$\fi}

\def\slch{{\cal SL}^c({\h})}
\def\slc{{\cal SL}^c}
\def\slqh{{\cal SL}({q,\h})}
\def\slcqh{{\cal SL}^c({q,\h})}
\def\swv{{\rm SW}(V_{(m|n)})}

\def\gh{\gg_{\h}}
\def\brpl{\{\,,\}_{PL}}
\def\brplr{\{\,,\}_{PL,r}}
\def\brr{\{\,,\}_{r}}
\def\brkks{\{\,,\}_{KKS}}
\def\brkksr{\{\,,\}_{KKS,r}}
\def\Tr{{\rm Tr}}
\def\lp{\Lambda_+(V)}
\def\lm{\Lambda_-(V)}
\def\lpm{\Lambda_{\pm}(V)}

\def\la{{\lambda}}
\def\al{{\alpha}}

\def\vl{{V_\lambda}}
\def\lb#1#2{{^#2\!#1}}

\def\be{\begin{equation}}
\def\ee{\end{equation}}

\begin{document}

\makeatletter
\renewcommand{\theequation}{{\thesection}.{\arabic{equation}}}
\@addtoreset{equation}{section} \makeatother

\title{Representation theory of (modified) Reflection Equation
Algebra of $GL(m|n)$ type}

\author{
\rule{0pt}{7mm} Dimitri
Gurevich\thanks{gurevich@univ-valenciennes.fr}\\
{\small\it USTV, Universit\'e de Valenciennes,
59304 Valenciennes, France}\\
\rule{0pt}{7mm} Pavel Pyatov\thanks{pyatov@thsun1.jinr.ru}\\
{\small\it Bogoliubov Laboratory of Theoretical Physics,
JINR, 141980 Dubna, Moscow region, Russia}\\
\rule{0pt}{7mm} Pavel Saponov\thanks{Pavel.Saponov@ihep.ru}\\
{\small\it Division of Theoretical Physics, IHEP, 142284
Protvino,
Russia} }

\maketitle

\begin{abstract}
Let $R:V^{\ot 2}\to V^{\ot 2}$ be a Hecke type solution of the
quantum Yang-Baxter equation (a Hecke symmetry). Then, the
Hilbert-Poincr\'e series of the associated $R$-exterior algebra
of
the space $V$ is a ratio of two polynomials of degree $m$
(numerator) and $n$ (denominator).

Assuming $R$ to be skew-invertible, we define a rigid quasitensor
category $\swv$ of vector spaces, generated by the space $V$ and
its
dual $V^*$, and compute certain numerical characteristics of its
objects. Besides, we introduce a braided bialgebra structure in
the
modified Reflection Equation Algebra, associated with $R$, and
equip
objects of the category $\swv$ with an action of this algebra. In
the case related to the quantum group $U_q(sl(m))$, we consider
the
Poisson counterpart of the modified Reflection Equation Algebra
and
compute the semiclassical term of the pairing, defined via the
categorical (or quantum) trace.

\end{abstract}

{\bf AMS Mathematics Subject Classification, 1991:} 17B37, 81R50


\section{Introduction}
\label{sec:int} Reflection Equation Algebra (REA) is a very
useful
tool of the theory of integrable systems with boundaries. It
derives
its name from an equation describing the factorized scattering on
a
half-line (cf. \cite{C}, where the REA depending on a spectral
parameter was first introduced).

Nowadays, different types of the REA are known (cf. \cite{KS}),
which have applications in mathematical physics and
non-commutative
geometry.

The REA related to the Drinfeld-Jimbo Quantum Group (QG) $\Uq$
appears in constructing a $q$-analog of differential calculus on
the
groups $GL(m)$ and $SL(m)$, where it was treated to be a
$q$-analog
of the exponential of vector fields (cf. \cite{FP}).

In the case related to the QG $U_q(\frak g)$, an appropriate
quotient of the REA can be treated as a deformation of the
coordinate ring $\K[G]$ where $G$ is the Lie group, corresponding
to
a classical Lie algebra $\frak g$. The Poisson bracket
corresponding
to this deformation was introduced by
M.Semenov-Tian-Shansky\footnote{Note that on any classical Lie
group
$G$ there exists another  Poisson bracket due to E.Sklyanin. Its
quantum analog is an appropriate quotient of the so-called RTT
algebra (cf. \cite{FRT}). These two quantum analogs of the space
$\K[G]$ are related by a transmutation procedure introduced by S.
Majid (cf. \cite{M} and references therein). Nowadays, there
exists
their universal treatment based on pairs of so-called compatible
braidings (cf. \cite{IOP,GPS1,GPS2}).}.

Though the best known REA is related to the QG $U_q(\frak g)$,
such
an algebra can be associated to any braiding $R:V^{\otimes 2}\to
V^{\otimes 2}$, where $V$ is a finite dimensional linear space
over
the ground field\footnote{Mainly we are dealing with $\K=\C$ but
sometimes $\K=\R$ is allowed.} $\Bbb K$ and $R$ is an invertible
solution of the quantum Yang-Baxter equation \be
R_{12}R_{23}R_{12}=R_{23}R_{12}R_{23}\,. \label{YB} \ee Here the
indices of $R$ relate to the space (or spaces) in which the
operator
is applied. Thus, $R_{12}$ and $ R_{23}$ are the following
operators
in the space $V^{\ot 3}$: $R_{12}=R\ot I$, $R_{23}=I \ot R$.

In the present paper we deal with {\em Hecke type} solutions of
the
Yang-Baxter equation (\ref{YB}) which satisfy the following
condition \be (R-q\, I)(R+q^{-1}\,I)=0\,, \label{Hec} \ee where
the
nonzero parameter $q\in {\Bbb K}$ is assumed to be {\em generic}.
By
definition, this means, that the values of $q$ do not belong to a
countable set of the roots of unity: $q^k\not=1$, $k=2,3,...$
(whereas the value $q=1$ is not excluded). Consequently,
$$
k_q:={q^k-q^{-k}\over q-q^{-1}}\not=0,\qquad \forall\,k\in {\Bbb
N}\,,
$$
$k_q$ being a a $q$-analog an integer $k$. In what follows, a
braiding satisfying relation (\ref{Hec}) will be called {\em a
Hecke
symmetry}.

Especially, we are interested in families of Hecke symmetries
$R_q$
analytically depending on the parameter $q$ in a neighbourhood of
$1\in {\Bbb K}$ in such a way, that for $q=1$ the symmetry
$R=R_1$
is involutive: $R^2=I$.

The well known example of such a family is the $U_q(sl(m))$
Drinfeld-Jimbo braidings \be R_q=\sum_{i,j =
1}^mq^{\delta_{ij}}\,h_i^j\ot h_j^i+\sum_{i<j}^m(q-q^{-1})\,
h_i^i\ot h_j^j \label{Rmat} \ee where the elements $h_i^j$ form
the
natural basis in the space of left endomorphisms of $V$, that is
$h_i^j(x_k)=\delta^j_k\,x_i$ in a fixed basis $\{x_k\}$ of the
space
$V$. Note that for $q=1$ the above braiding $R$ equals the usual
flip $P$.

The Hecke symmetry (\ref{Rmat}) and all related objects will be
called {\em standard}. However, a large number of Hecke
symmetries
different from the standard one are known, even those which are
not
deformations of the usual flip (cf. \cite{G3}).

Let us consider the REA corresponding to the standard
$U_q(sl(m))$
Hecke symmetry (\ref{Rmat}) in more detail. This algebra
possesses
some very important properties, in contrast with the REA related
to
other quantum groups $U_q(\gg)$, $\gg\not=sl(m)$.

First of all, it is a $q$-deformation of the commutative algebra
$\Sym(gl(m))=\K[gl(m)^*]$ (so, we get a deformation algebra
without taking
any additional quotient). Second, by a linear shift of REA
generators (proportional to a parameter $\hbar$), we come to
quadratic-linear commutation relations for the shifted
generators.
In this basis the REA can be treated as a "double deformation" of
the initial commutative algebra $\K[gl(m)^*]$. We refer to this
form
of the REA as {\em modified Reflection Equation Algebra} (mREA)
and we
 denote it $\lrqh$. By specializing  $\h=0$ we return to the
(non-modified) REA $\lrq$.

The specialization of the algebra $\lrqh$ at $q=1$ gives the
enveloping algebra $U(gl(m)_\h)$ where the notation  $\gh$ means
that the bracket $[\,,\,]$ of a Lie algebra $\gg$ is replaced by
$\h[\,,\,]$. (Note, that this fact was observed in \cite{IP}.)
The
commutative algebra $\K[gl(m)^*]$ is obtained by the double
specialization of the algebra $\lrqh$ at $\h=0$ and $q=1$.

Being equipped with the $\Uq$-module structure, the algebra
$\lrqh$
(as well as $\lrq$) is $\Uq$-equivariant (or covariant). This
means
that
$$
M(x\cdot y)=M_{(1)}(x)\cdot M_{(2)}(y)\,,\quad\forall\, M\in
\Uq\,,\quad \forall\, x,y\in \lrqh
$$
where we use the Sweedler's notation for the quantum group
coproduct
$\Delta(M)=M_{(1)}\ot M_{(2)}$.

The Poisson counterpart of the above double deformation of the
algebra $\K[gl(m)^*]$ is the Poisson pencil \be
\brplr=a\,\brpl+b\,\brr\,,\quad a,b\in \K \label{pp} \ee where
$\brpl$ is the linear Poisson-Lie bracket related to the Lie
algebra
$gl(m)$ and $\brr$ is a natural extension  of the
Semenov-Tian-Shansky bracket  on the linear space $gl(m)^*$.

We consider these Poisson structures and briefly discuss their
role
in defining a "quantum orbit" ${\cal O} \subset gl(m)^*$ in
Section
\ref{sec:tr-qua}. Taking a two-dimensional sphere as an example,
we suggest a method of constructing such quantum orbits. In
contrast
with other definitions of quantum homogeneous spaces, our quantum
orbits are some quotients of the algebra $\lrqh$. They look like
the
"fuzzy sphere"
$$
\slch=U(su(2)_\h)/\langle C-c \rangle
$$
where $C$ is the quadratic Casimir element. As is known, there
exists a discrete series of numbers $c_k\in {\Bbb K}$ such that
any
algebra ${\cal SL}^{c_k}(\h)$ has a finite dimensional
representation in a linear space $V_k$ and the corresponding map
${\cal SL}^{c_k}(\h) \to \End(V_k)$ is an $su(2)$-morphism.

A similar statement is valid for the aforementioned quotients of
the
algebra $\lrqh$. However, the corresponding spaces $V_k$ become
objects of a quasitensor category. In such a category, an object
is
characterized by its quantum dimension which is defined via the
categorical (quantum) trace. A deformation of the usual trace is
one
of the main features of our approach to the quantum homogeneous
spaces. In Section \ref{sec:tr-qua} we describe the semiclassical
term of the paring defined via the quantum trace in the case of
the
standard Hecke symmetry.

In a similar way we treat other quasitensor categories generated
by
skew-invertible Hecke symmetries. Roughly speaking, we are
dealing
with three problems in the present paper. The first problem is
the
classification of all (skew-invertible) Hecke symmetries $R$. One
of
the main tools for studying this problem is the
Hilbert-Poincar\'e
(HP) series $P_-(t)$ corresponding to the "$R$-exterior algebra"
of
the space $V$ (its definition is presented in Section
\ref{sec:g-fo}). Though a classification of all possible forms of
the HP series $P_-(t)$ has not been found yet, it is known that
the
HP series $P_-(t)$ of any Hecke symmetry is a rational
function\footnote{The HP series corresponding to a
skew-invertible
Hecke symmetry is described in Section \ref{sec:g-fo}. When
$P_-(t)$
is a polynomial (in this case we say that $R$ is even) it can
drastically differ from the classical one $(1+t)^n$, $n=\dim V$.
Thus, in \cite{G3} all skew-invertible Hecke symmetries with
$P_-(t)=1+nt+t^2$ were classified. Besides, suggested in
\cite{G3}
was a way of "gluing" such symmetries which gives rise to
skew-invertible Hecke symmetries with other non-standard HP
series.}
\cite{H,D}. The ordered pair of integers $(m|n)$, where $m$
(resp.,
$n$) is the degree of the numerator $N(t)$ (resp., denominator
$D(t)$) of $P_-(t)$, plays an important role in the sequel and
will
be called the {\em bi-rank} of the Hecke symmetry $R$ (or of the
corresponding space $V$). This pair enters our notation of the
quasitensor Schur-Weyl category $\swv$ generated by $V$.

Constructing the category $\swv$ is the second problem we are
dealing with in this paper. The  objects of the category are
direct
sums of vector spaces $V_\la\ot V^*_\mu$. Here $V$ is the basic
vector space equipped with a skew-invertible Hecke symmetry $R$,
$V^*$ is its dual, and $\la$ and $\mu$ stand for arbitrary
partitions (Young diagrams) of positive integers. The map $V\to
\vl$
is nothing but the Schur functor corresponding to the Hecke
symmetry
$R$ (for its classical version cf. \cite{FH}). The map $V^* \to
V^*_\mu$ can be defined in a similar way. Note, that $\swv$ is a
monoidal quasitensor rigid category (as defined in \cite{CP}) but
it
is not abelian.

We compute some numerical characteristics of objects of this
category. Namely, we are interested in their dimensions
(classical
and quantum). In contrast with the classical dimensions, which
essentially depend on a concrete form of the initial Hecke
symmetry
and is expressed via the roots of above polynomials $N(t)$ and
$D(t)$, the quantum dimensions depend only on the bi-rank
$(m|n)$.
Moreover, in a sense, the category $\swv$ looks like the tensor
category of $U(gl(m|n))$-modules.

The third problem elaborated below is in constructing the
representations of the mREA $\lrqh$ in the category $\swv$. Since
for $q\not=1$ the algebra $\lrqh$ is isomorphic to the
non-modified
REA $\lrq$ (in fact, we have the same algebra written in two
different bases), we automatically get a representation category
of the latter algebra\footnote{Since for $q=1$ the isomorphism
$\lrqh\cong \lrq$ breaks, we prefer to consider these algebras
separately and use different names for them.}. Note, that certain
representations of the REA have already been known, mainly for
the even
case (the bi-rank $(m|0)$) \cite{K,Mu1,GS2,S}. In contrast with
those papers, we here consider the mREA $\lrqh$ connected with a
general type skew-invertible Hecke symmetry $R$ of the bi-rank
$(m|n)$ and equip objects of the category $\swv$ with the
$\lrqh$-module structure. Note, that all the corresponding
representations are equivariant (see Section \ref{sec:7}).

A particular example we are interested in is the "adjoint"
representation. By this we mean a representation $\rho_{ad}$ of
the
mREA $\lrqh$ in the linear span of its generators. In the case,
when
a Hecke symmetry is a super-flip in a $Z_2$-graded linear space V
$$
R:\;\vv\to\vv\,,\qquad R(x\ot
y)=(-1)^{\overline{x}\,\overline{y}}y\ot x\,,
$$
where $x$ and $y$ are homogenous elements of $V$ and $\bar z$
denotes the parity (grading) of a homogeneous element $z$, the
mREA
becomes the enveloping algebra $U(gl(m|n))$ and the
representation
$\rho_{ad}$ coincides with the usual adjoint one. This is one of
the
reasons why we treat the mREA $\lrqh$ as a suitable analog of the
enveloping algebra. Moreover, in the case of involutive
skew-invertible Hecke symmetry, the corresponding mREA becomes
the
enveloping algebra of a generalized Lie algebra $\End(V)$ as is
explained in Section \ref{sec:mrea}. Such  algebras were
introduced in \cite{G1}.

The other property that makes the mREA  similar to the
enveloping algebra of a generalized Lie algebra (in particular, a
super Lie algebra) is its braided bialgebra structure. Such a
structure is determined by a coproduct $\Delta$ and a counit
$\varepsilon$. On the generators of mREA (organized into a matrix
$L$ (see Section \ref{sec:7})) the coproduct reads
$$
\Delta(L)=L\ot 1+1\ot L-(q-q^{-1})L\ot L
$$
and coincides with the coproduct of the enveloping algebra of the
(generalized) Lie algebra at $q=1$. Note, that though we do not
define an antipode in the algebra $\lrqh$, the category $\swv$ of
its representations is closed.

In addition to the $\lrqh$-module structure, the objects of the
Schur-Weyl category, corresponding to the standard Hecke symmetry
(\ref{Rmat}), can be equipped with the action of the QG $\Uq$.
Besides, the $q$-analogs of super-groups (cf. \cite{KT}) can also
be
represented in the corresponding Schur-Weyl category. (Suggested
in
\cite{Z} is another way of constructing the representations of
$q$-deformed algebras $U(gl(m|n)$ which is based on the
triangular
decomposition.) Nevertheless, in general we know no explicit
construction of the QG type algebra for a skew-invertible Hecke
symmetry\footnote{An attempt of explicit description of such an
object for some even non-quasiclassical Hecke symmetries was
undertaken in \cite{AG}.} whereas  the mREA can be defined for
any
skew-invertible Hecke symmetry.

The mREA has one more advantage compared with the QG or their
super-analogs. It is a more convenient tool for the explicit
construction of projective modules over quantum orbits in the
frameworks of approach suggested in \cite{GS1,GS3}. We plan to
turn
to these objects in a general (not necessarily even) case in our
subsequent publications.

To complete the Introduction, we would like to emphasize a
difference between the Hecke type braidings and the
Birman-Murakami-Wenzl ones (in particular, those coming from the
QG
of $B_n$, $C_n$ and $D_n$ series). In the latter case it is not
difficult to define a "braided Lie bracket" in the space
$\End(V)$
(cf. \cite{DGG}) and introduce the corresponding "enveloping
algebra". But this "enveloping algebra" is not a deformation of
its
classical counterpart and therefore is not an interesting object
from our viewpoint.

The paper is organized as follows. In the next Section we
reproduce
some elements of $R$-technique  which form the base of subsequent
computations of some interesting numerical characteristics of
objects involved (the most cumbersome part of the computations is
placed in Appendix). Section \ref{sec:g-fo} is devoted to the
classification of (skew-invertible) Hecke symmetries. In Section
\ref{sec:SWH} we construct the Schur-Weyl category $\swv$
generated
by the space $V$. Our main object, the mREA $\lrqh$, is
introduced
in Section \ref{sec:mrea} where we also study its deformation
properties. In Section \ref{sec:7} we equip the mREA with a
braided
bialgebra structure which allows us to define an equivariant
action
of the algebra $\lrqh$ on each object of the category $\swv$.
There
we also present our viewpoint on definition of braided (quantum)
Lie
algebras. Section \ref{sec:tr-qua} is devoted to study of some
semiclassical structures.
\bigskip

\noindent{\bf Acknowledgement.} We would like to thank the
Max-Planck-Institut f{\"u}r Mathematik, where this work was
written,
for the warm hospitality and stimulating atmosphere. The work of
D.G. was partially supported by the grant ANR-05-BLAN-0029-01,
the
work of P.P. and P.S. was partially supported by the RFBR grant
05-01-01086.

\section{Elements of $R$-technique}
\label{sec:2}

By $R$-technique we mean  computational methods based on general
properties of braidings (in particular, Hecke symmetries)
regardless
of their concrete form. We are mostly interested in the so-called
skew-invertible braidings since they enable us to define
numerical
characteristics of Hecke symmetries and related objects.

A braiding $R$ (see (\ref{YB})) is called {\em skew-invertible}
if
there exists an endomorphism $\Psi: V^{\otimes 2}\rightarrow
V^{\otimes 2}$ such that \be \Tr_{(2)} \,
R_{12}\,\Psi_{23}=P_{13}=\Tr_{(2)} \, \Psi_{12}\,R_{23}
\label{psi}
\ee where the symbol $\Tr_{(2)}$ means calculating trace in the
second factor of the tensor product $V^{\otimes 3}$. Hereafter
$P$
stands for the usual flip $P(x\ot y)=y\ot x$.

Fixing bases $\{x_i\}$ and $\{x_i\ot x_j\}$ in $V$ and
$V^{\otimes
2}$ respectively, we identify $R$ (resp., $\Psi$) with a matrix
$\|R_{\;ij}^{kl}\|$ (resp., $\|\Psi_{\;ij}^{kl}\|$): \be R
(x_{i}\ot
x_{j})=x_{k}\ot x_{l} \, R_{\;ij}^{kl} \label{flip} \ee where the
upper indices mark the rows of the matrix and from now on the
summation over the repeated indices is understood.

Being written in terms of matrices, relation (\ref{psi}) reads
$$
R_{\;jb}^{ia}\,\Psi_{\;ak}^{bl}=\delta_k^i \delta_j^l=
\Psi_{\;jb}^{ia}\,R_{\;ak}^{bl}\,.
$$

Using $\Psi$ we define two endomorphisms $B$ and $C$
of the space $V$
$$
B(x_i)= x_j B^j_{\;i},\qquad C(x_i)= x_j C^j_{\;i},
$$
where \be B_{\;i}^j:=\Psi^{kj}_{\;ki},\qquad C_{\;i}^j:=
\Psi^{jk}_{\;ik}, \label{def:B-C} \ee that is
$$
B:=\Tr_{(1)}\,\Psi,\qquad C:=\Tr_{(2)}\,\Psi.
$$

If the operator $B$ (or $C$) is invertible, then the
corresponding
braiding $R$ is called {\it strictly} skew-invertible. As was
shown in \cite{O}, $R$ is strictly skew-invertible iff $R^{-1}$
is
skew-invertible and, besides, the invertibility of $B$ leads to
the
invertibility of $C$ and vice versa.

A well known important example of a strictly skew-invertible
braiding is the super-flip $R$ on a super-space $V=V_0\oplus
V_1$,
where $V_0$ and $V_1$ are respectively the even and odd
components
of $V$. In this case the operators $B$ and $C$ are called {\em
the
parity operators} and their explicit form is as follows
$$
B(z) = C( z) = z_0-z_1,\quad \forall\, z\in V,
$$
where $z_0(z_1)$ is the even (odd) component of $z=z_0+z_1$.

Let $R$ be a skew-invertible braiding. Listed below are some
useful
properties of the corresponding endomorphisms $\Psi$, $B$ and
$C$.

\begin{enumerate}
\item
\parbox[b]{0.9\hsize}{\centerline{$\Tr\,B=\Tr\, C,$}}
\par
\be \qquad \Tr_{(2)}\,B_{2}\,R_{21}=\Tr_{(2)}\,C_{2}\,
R_{12}=I_1,
\label{trR} \ee where $I$ is the identical automorphism of $V$.
These relations directly follow from definitions (\ref{psi}) and
(\ref{def:B-C}).

\item The endomorphisms $B$ and $C$ commute and their product
is a scalar operator \be B\, C=C\, B = \nu\, I\label{BC}, \ee
where
the numeric factor $\nu$ is nonzero iff the braiding $R$ is
strictly
skew-invertible (in particular, if $R$ is a skew-invertible Hecke
symmetry).

\item
The matrix elements of $B$ and $C$ realize a one-dimensional
representation of the so-called $RTT$ algebra, associated with
$R$
(cf. \cite{FRT}), that is \be
R_{12}B_{1}B_{2}=B_{1}B_{2}R_{12},\qquad
R_{12}C_{1}C_{2}=C_{1}C_{2}R_{12}\,. \label{bc-rtt} \ee As a
direct
consequence of the above relations, we have
$$
\Tr_{(12)}(B_{1}B_{2}R_{12}
X_{12}R_{12}^{-1})=\Tr_{(12)}(B_{1}B_{2}R_{12}^{-1}
X_{12}R_{12})=\Tr_{(12)}(B_{1}B_{2}X_{12}),
$$
$$
\Tr_{(12)}(C_{1}C_{2}R_{12}X_{12}R_{12}^{-1})=
\Tr_{(12)}(C_{1}C_{2}R_{12}^{-1}X_{12}R_{12})=
\Tr_{(12)}(C_{1}C_{2}X_{12})
$$
where $X\in{\rm End}(V^{\otimes 2})$ is an arbitrary endomorphism
and $\Tr_{(12)}(\dots)=\Tr_{(1)}(\Tr_{(2)}(\dots))$.

\item The following important relations were proved in
\cite{OP2,S} \be
\begin{array}{lcl}
B_{1}\Psi_{12}=R_{21}^{-1}B_{2},&\qquad&
\Psi_{12}B_{1}=B_{2}R_{21}^{-1},\\
\rule{0pt}{5mm} C_{2}\Psi_{12}=R_{21}^{-1}C_{1},&\qquad&
\Psi_{12}C_{2}=C_{1}\,R_{21}^{-1},
\end{array}
\label{dva} \ee where $R_{21}=PR_{12}P$. In case $\nu\not=0$,
only
one of the lines above is independent due to relation (\ref{BC}).

Therefore, for an arbitrary endomorphism $X\in {\rm End}(V)$ we
obtain
\begin{eqnarray}
&&\Tr_{(1)}(B_{1}R_{12}X_{2}R_{12}^{-1})=
\Tr_{(1)}(B_{1}R_{12}^{-1}X_{2}R_{12})=
\Tr(B X)\, I_2,\nonumber\\
&&\Tr_{(2)}(C_{2}R_{12}X_{1}R_{12}^{-1})=
\Tr_{(2)}(C_{2}R_{12}^{-1}X_{1}R_{12})= \Tr(C
X)\,I_1.\label{trRC}
\end{eqnarray}
\end{enumerate}

This completes the list of technical facts to be used in the text
below.

\section{The general form of a Hecke symmetry}
\label{sec:g-fo}

In this Section we study the classification problem of
(skew-invertible) Hecke symmetries. Our presentation is based on
the
theory of the $A_{k-1}$ series Hecke algebras and their
$R$-matrix
representations. As a review of the subject we can recommend the
work \cite{OP1}. Some necessary facts of the mentioned theory are
given in Appendix for the reader's convenience.

Given a Hecke symmetry $R:V^{\otimes 2}\rightarrow V^{\otimes
2}$,
we consider the $R$-symmetric $\lp$ and the $R$-skew-symmetric
$\lm$
algebras of the space $V$, which by definition are the following
quotients \be \lpm:=T(V)/\langle({\rm Im}(q^{\pm 1}\,I_{12}\mp
R_{12})\rangle, \qquad I_{12} = I\otimes I. \label{La-pm} \ee
Hereafter $T(V)$ stands for the free tensor algebra of the space
$V$
and $\langle J\rangle$ denotes the two-sided ideal generated in
this
algebra by a subset $J\subset T(V)$.

Then, we consider the Hilbert-Poincar\'e (HP) series of the
algebras
$\lpm$ \be P_{\pm}(t):=\sum_{k\ge 0}
t^k\,\dim\,\Lambda^k_{\pm}(V),
\label{hp-ser} \ee where $\Lambda^k_{\pm}(V)\subset
\Lambda_{\pm}(V)$ is the homogenous component of degree $k$.

The following proposition plays a decisive role in the
classification of all possible forms of the Hecke symmetries.
\begin{proposition}
\label{pro:1} Consider an arbitrary Hecke symmetry $R$,
satisfying
{\rm (\ref{YB})} and {\rm (\ref{Hec})} at a generic value of the
pa\-ra\-me\-ter $q$. Then the following properties hold true.
\begin{enumerate}
\item[\rm 1.]
The HP series $P_{\pm}(t)$ obey the relation
$$
P_+(t)\,P_-(-t)=1.
$$

\item[\rm 2.]
The HP series $P_-(t)$ {\rm ({\it and hence} $P_+(t)$)} is a
rational function of the form: {\rm \be P_-(t)=\frac{N(t)}{D(t)}=
\frac{1+a_1\,t+...+a_m\,t^m}{1-b_1\,t+...+(-1)^n\, b_n\,t^n} =
\frac{\prod_{i=1}^m(1+x_it)}{\prod_{j=1}^n(1-y_jt)}\,, \label{p-}
\ee} where the coefficients $a_i$ and $b_i$ are positive
integers,
the polynomials $N(t)$ and $D(t)$ are mutually prime, and all
real
numbers $x_i$ and $y_i$ are positive.

\item[\rm 3.]
If, in addition, the Hecke symmetry is skew-invertible, then the
polynomials $N(t)$ and $D(-t)$ are reciprocal\footnote{Recall,
that
a polynomial $p(t) = c_0 + c_1t +\dots +c_nt^n$ with real
coefficients $c_i$ is called {\it reciprocal} if $p(t) =
t^np(t^{-1})$ or, equivalently, $c_i = c_{n-i}$, $0\le i\le n$.}.
\end{enumerate}
\end{proposition}

The first item of the above list was proved in \cite{G2}, the
second
and the third ones --- in \cite{H,Da} and \cite{DH}.

\begin{definition}
\label{def:2} {\rm Let $R:V^{\otimes 2}\rightarrow V^{\otimes 2}$
be
a skew-invertible Hecke symmetry and let $m$ (resp., $n$) be the
degree of the numerator $N(t)$ (resp., the denominator $D(t)$) of
the
HP series $P_-(t)$. The ordered pair of integers $(m|n)$ will be
called {\it the bi-rank} of $R$. If $n=0$ (resp., $m=0$), the
Hecke
symmetry will be called {\it even} (resp., {\it odd}). Otherwise
we
say that $R$ is of the general type.}
\end{definition}

\begin{remark}\rm
In the sense of the above definition, any skew-invertible Hecke
symmetry is a generalization of the super-flip for which
$P_-(t)=(1+t)^m(1-t)^{-n}$, where $m=\dim\,V_0$, $n=\dim\,V_1$.
Such
a treatment of Hecke symmetries is also motivated by similarity
of
the corresponding Schur-Weyl categories (see below).
\end{remark}

Now we obtain some important consequences of Proposition
\ref{pro:1}. Let $R$ be a Hecke symmetry of the bi-rank $(m|n)$.
As
is known, the Hecke symmetry $R$ allows to define a
representations $\rho_R$ of the $A_{k-1}$ series Hecke algebras
$H_k(q)$, $k\ge 2$, in homogeneous components  $V^{\otimes p}
\subset T(V)$, $\forall\, p\ge k$
$$
\rho_R:\; H_k(q)\rightarrow {\rm End}(V^{\otimes p}),\quad p\ge
k\,.
$$
Explicitly, these representations are given in (\ref{r-rep-H}) of
Appendix.

Under the presentation $\rho_R$, the primitive idempotents
$e_a^\lambda\in H_k(q)$, $\lambda\vdash k$, convert to the
projection operators \be E^{\lambda}_a(R)=\rho_R(e^\lambda_a)\in
\End(V^{\otimes p}), \quad p\ge k\,, \label{Y-proj} \ee where the
index $a$ enumerates the standard Young tableaux $(\lambda,a)$,
which can be  constructed for a given partition $\lambda\vdash
k$.
The total number of the standard Young tableaux corresponding to
the
partition $\lambda$ is denoted as $d_\lambda$.

Under the action of these projectors the spaces $V^{\otimes p}$,
$p\ge 2$, are expanded into the direct sum \be V^{\otimes p} =
\bigoplus_{\lambda\vdash p} \bigoplus_{a=1}^{d_\lambda}
V_{(\lambda,a)}, \qquad V_{(\lambda,a)} = {\rm Im}(E_a^\lambda).
\label{V-decom} \ee Due to relation (\ref{id-conn}), the
projectors
$E^\lambda_a$ with different $a$ are connected by invertible
transformations and, therefore, all spaces $V_{(\lambda,a)}$ with
fixed $\lambda$ and different $a$ are isomorphic.

At a generic value of $q$, the Hecke algebra $H_k(q)$ is known to
be
isomorphic to the group algebra ${\Bbb K}[{\frak S}_k]$
\cite{We}.
Basing on this fact, we can prove the following result
\cite{GLS1,H}
\be V_{(\lambda,a)}\otimes V_{(\mu,b)} = \bigoplus_{\nu}
\bigoplus_{d_{ab}\in I_{ab}} V_{(\nu,d_{ab})}\cong
\bigoplus_{\nu}
c^\nu_{\lambda\mu}\, V_{(\nu,d_0)}\,,\qquad \lambda\vdash p,
\,\mu\vdash k\,, \nu\vdash(p+k)\,, \label{v-lit-r} \ee where the
integers $c^\nu_{\lambda\mu}$ are the Littlwood-Richardson
coefficients, the tableau index $d_{ab}$ takes the values form a
subset $I_{ab}\subset \{1,2,\dots,d_\nu\}$, which depends on the
values of the indices $a$ and $b$. The number $d_0$ in the last
equality stands for the index of an arbitrary fixed tableau from
the
set $(\nu,d)$, $1\le d\le d_\nu$. This equality has the following
meaning. Though the summands $ V_{(\nu,d_{ab})}$ do depend on the
values of $a$ and $b$, the {\it total} number of these summands
(the
cardinality of $I_{ab}$) depends only on the partitions
$\lambda$,
$\mu$ and $\nu$ and is equal to the Littlewood-Richardson
coefficient $c^\nu_{\lambda\mu}$. Therefore, due to isomorphism
$V_{(\nu,d_{ab})} \cong V_{(\nu,d_0)}$, we can replace the sum
over
$d_{ab}$ by the space $V_{(\nu,d_0)}$ with the corresponding
multiplicity $c^\nu_{\lambda\mu}$ (cf. \cite{GLS1}).

A particular example of the spaces $V_{(\lambda,a)}$ is the
homogeneous components $\Lambda_+^k(V)$ and $\Lambda_-^k(V)$ of
the
algebras $\Lambda_\pm(V)$ (\ref{La-pm}). They are images of the
projectors $E^{(k)}$ and $E^{(1^k)}$, corresponding to one-row
and
one-column partitions $(k)$ and $(1^k)$ respectively. This
important
fact allows us to calculate the dimensions (over the ground field
$\Bbb K$) of all spaces $V_{(\lambda,a)}$, provided that the
Poincar\'e series $P_-(t)$ is known. Since all the spaces
$V_{(\lambda,a)}$ corresponding to the same partition $\lambda$
are
isomorphic, we denote their $\Bbb K$-dimensions by the symbol
$\dim
V_\lambda$.

In the sequel, the following corollary of Proposition \ref{pro:1}
will be useful.
\begin{corollary}
\label{cor:4-1} Let $R$ be a Hecke symmetry of the bi-rank
$(m|n)$,
the Poincar\'e series of $\Lambda_-(V)$ being given by {\rm
(\ref{p-})}. Then for the partitions $(k)$ and $(1^k)$, $k\in
{\Bbb
N}$, the dimensions of the spaces $V_{(k)}$ and $V_{(1^k)}$ is
determined by the formulae {\rm
\begin{eqnarray}
\dim V_{(k)} &=& s_{(k)}(x|y):=\sum_{i=0}^kh_{i}(x)e_{k-i}(y),
\label{dim-row}\\
\dim V_{(1^k)} &=&
s_{(1^k)}(x|y):=\sum_{i=0}^ke_{i}(x)h_{k-i}(y),
\label{dim-col}
\end{eqnarray}}
where $h_i$ and $e_i$ are respectively the complete symmetric and
elementary symmetric functions of their arguments.
\end{corollary}

\noindent {\bf Proof.\ } We prove only the first of the above
formulae since the second one can be proved in the same way.
Since
$V_{(k)} = \Lambda_+^k(V)$, the dimension of $V_{(k)}$ can be
found
as an appropriate derivative of the Poincar\'e series $P_+(t)$
$$
\dim V_{(k)} = \frac {1}{k!}\,\frac{d^k}{dt^k} P_+(t)_{|_{t=0}}.
$$
Using $P_+(t)P_-(-t) = 1$ (see Proposition \ref{pro:1})
and relation (\ref{p-}) we  present $P_+(t)$ in the form
$$
P_+(t) = \prod_{i=1}^n(1+y_it)\prod_{j=1}^m\frac{1}{(1-x_jt)} =
{\cal E}(y|t){\cal H}(x|t),
$$
where ${\cal E}(\cdot)$ and ${\cal H}(\cdot)$ stands for the
generating functions of the elementary and complete symmetric
functions in the finite set of variables \cite{Mac}:
\begin{eqnarray*}
&&e_k(y) = \sum_{1\le i_1<\dots<i_k\le n}y_{i_1}\dots y_{i_k} =
\frac {1}{k!}\,\frac{d^k}{dt^k} {\cal E}(y|t)_{|_{t=0}}\\
&&h_k(x) = \sum_{1\le j_1\le\dots\le j_k\le m}x_{j_1}\dots
x_{j_k} =
\frac {1}{k!}\,\frac{d^k}{dt^k} {\cal H}(x|t)_{|_{t=0}}\,.
\end{eqnarray*}
Calculating the $k$-th derivative of $P_+(t)$ at $t=0$ we get
 (\ref{dim-row}).\hfill\rule{6.5pt}{6.5pt}
\medskip

Note, that polynomials $s_{(k)}(x|y)$ and $s_{(1^k)}(x|y)$
defined
in (\ref{dim-row}) and (\ref{dim-col}) belong to the class of
super-symmetric polynomials in $\{x_i\}$ and $\{y_j\}$. By
definition \cite{St}, a polynomial $p(u|v)$ in two sets of
variables
is called {\it super-symmetric} if it is symmetric with respect
to
any permutation of arguments $\{u_i\}$ as well as of arguments
$\{v_j\}$ and, additionally, on setting $u_1 = v_1 = t$ in
$p(u|v)$
one gets the result independent of $t$. Evidently, the
polynomials
in question satisfy this definition if we set, for example, $u =
x$,
$v = -y$.

Actually, the set of polynomials $s_{(k)}(x|y)$ (respectively
$s_{(1^k)}(x|y)$),  $k\in {\Bbb N}$,  are super-symmetric analogs
of
complete symmetric (respectively elementary symmetric) functions
in
finite numbers of variables. In particular, they generate the
whole
ring of super-symmetric polynomials in variables $\{x_i\}$ and
$\{y_j\}$. The $\Bbb Z$-basis of this ring is formed by the Schur
super-symmetric functions $s_\lambda(x|y)$ which can be expressed
in
terms of $s_{(k)}$ (or $s_{(1^k)}$) through Jacobi-Trudi
relations
\cite{Mac}. The Schur super-symmetric functions determine the
value
of dimensions $\dim V_{\lambda}$. In order to formulate the
corresponding result we need one more definition.
\begin{definition}\rm
\label{def:hook} {\bf (\cite{BR})} Given two arbitrary integers
$m\ge 0$ and $n\ge 0$, consider a partition $\lambda =
(\lambda_1,\lambda_2,\dots)$, satisfying the following
restriction
$\lambda_{m+1}\le n$. The (infinite) set of all such partitions
are
denoted as ${\sf H}(m,n)$ and any partition $\lambda\in {\sf
H}(m,n)$ will be called a hook partition of the type ${\sf
H}(m,n)$.
\end{definition}

\begin{proposition}{\bf (\cite{H})}
\label{pro:dim-Vl} Let $R$ be a Hecke symmetry of the bi-rank
$(m|n)$. Then the dimensions $\dim V_\lambda$ of spaces in
decomposition {\rm (\ref{V-decom})} are determined by the rules:
\begin{enumerate}
\item[{\rm 1.}] For any $\lambda=(\lambda_1,\dots ,\lambda_k)\in
{\sf H}(m,n)$ the dimension $\dim V_\lambda\not = 0$ and is given
by
the formula {\rm \be \dim V_\lambda = s_\lambda(x|y)\,.
\label{dim-Vl} \ee} Here
$$
s_\lambda (x|y)= \det\| s_{(\lambda_i-i+j)}(x|y) \|_{1\le i,j\le
k}\,.
$$
where $s_{(k)}(x|y)$ is defined in {\rm (\ref{dim-row})} for
$k\ge
0$ and $s_{(k)}:=0$ for $k<0$.
\item[{\rm 2.}] For arbitrary partition $\lambda$ we have
$$
\dim V_\lambda = 0\; \Leftrightarrow \;\lambda\not\in {\sf
H}(m,n).
$$
\end{enumerate}
\end{proposition}

\noindent {\bf Proof.} Taking into account that
$$
\dim (U\otimes W) = \dim U\dim W, \qquad \dim(U\oplus W)= \dim U
+
\dim W
$$
and calculating dimensions of the spaces in the both sides of
(\ref{v-lit-r}) we find
$$
\dim V_\lambda\dim V_\mu = \sum_{\nu}c^\nu_{\lambda\mu} \dim
V_\nu.
$$
Now the result (\ref{dim-Vl}) is a direct consequence of an
inductive procedure based on Corollary \ref{cor:4-1} (cf., for
example, \cite{GPS2}).

The second claim can be deduced from the properties of the Schur
functions $s_\lambda(x|y)$ established in \cite{BR} (also cf.
\cite{H}).\hfill \rule{6.5pt}{6.5pt}
\medskip

To end the Section, we present one more important numerical
characteristic of the Hecke symmetry which can be expressed in
terms
of its bi-rank.
\begin{proposition}
\label{teo:3} Let $R$ be a skew-invertible Hecke symmetry with
the
bi-rank $(m|n)$. Then {\rm \be Tr\, B=\Tr\, C= q^{n-m}(m-n)_q.
\label{trBC} \ee }
\end{proposition}
The proof of the theorem is rather technical and is placed in the
Appendix.

\begin{corollary}
\label{cor:4} For a skew-invertible Hecke symmetry with the
bi-rank
$(m|n)$, the factor $\nu$ in {\rm (\ref{BC})} equals $q^{2(n-m)}$
that is
$$
BC = CB = q^{2(n-m)}I\,.
$$
\end{corollary}
{\bf Proof.\ } First, observe, that if  $R$ is a skew-invertible
Hecke symmetry, then the same is true for the operator
$R_{21}=PR_{12}P$, and therefore
$$
R_{21}^{-1}= R_{21}-(q-q^{-1})\,I_{21}.
$$
Applying $\Tr_{(2)}$ to the first formula in (\ref{dva}), we have
$$
\begin{array}[b]{lcl}
B_{1}C_{1} \hspace*{-5pt}&=&\hspace*{-5pt}
\Tr_{(2)}\,(B_{1}\Psi_{12})=\Tr_{(2)}(R_{21}^{-1}B_{2})\\
\rule{0pt}{6mm} &=&\hspace*{-5pt}
\Tr_{(2)}((R_{21}-(q-q^{-1})\,I_{21})B_{2}) =
I_{1}-(q-q^{-1})\,I_{1}\Tr(B) =q^{2(n-m)}\, I_{1}.
\end{array}
\eqno{\rule{6.5pt}{6.5pt}}
$$

\section{The quasitensor category $\swv$}
\label{sec:SWH}

Our next goal is to construct the quasitensor {\it Schur-Weyl
category} $\swv$ of vector spaces, generated by the space $V$
equipped with a skew-invertible Hecke symmetry $R$ of the bi-rank
$(m|n)$. The objects of this category possess the module
structure
over the {\it reflection equation algebra}, which will be
considered
in detail in the next Sections.

In constructing the above mentioned category we proceed
analogously to the paper \cite{GLS1}, where such a category was
constructed for an even Hecke symmetry of the bi-rank $(m|0)$. A
peculiarity of the even case is that the space $V^*$, dual to
$V$,
can be identified with a specific object $V_{(1^{m-1})}$ (see
(\ref{V-decom}) for the definition of $\vl$) of the category.
This
property ensures the category constructed in \cite{GLS1} to be
rigid\footnote{Recall, that a (quasi)tensor category of vector
spaces is {\it rigid} if to any of its objects $U$ there
corresponds
a dual object $U^*$ such that the maps $U\otimes U^*\rightarrow
{\Bbb K}$  and $U^*\otimes U\rightarrow {\Bbb K}$ are categorical
morphisms.}.

It is not so in the case of general bi-rank $(m|n)$, and we have
to
properly enlarge the category by adding the dual spaces to all
objects. This requires, in turn, a consistent extending of the
categorical braidings to the dual objects and defining  the
invariant pairings. In the present Section we elaborate these
problems in detail.

So, let $R$ be a skew-invertible Hecke symmetry of bi-rank
$(m|n)$, which, upon fixing a basis $\{x_i\}_{1\le i\le N}$ of
the
space $V$, $\dim V = N$, is represented by the matrix
(\ref{flip}).
Also, we introduce the dual vector space $V^*$ and choose the
basis
$\{x^i\}_{1\le i\le N}$ in $V^*$ dual to $\{x_i\}$ with respect
to
the nondegenerate bilinear form \be \langle\,,\rangle_r:\;
V\otimes
V^*\rightarrow {\Bbb K}\,, \qquad \langle x_i,x^j\rangle_r =
\delta^j_i\,. \label{A1} \ee The subscript $r$ (right) refers to
the
order of arguments in the form $\langle\,,\,\rangle_r$: the
vectors
of the dual space $V^*$ stand on the right of the vectors of $V$.

By definition, the dual space to the tensor product $U\otimes W$
is
$W^*\otimes U^*$:
$$
\langle U\otimes W, W^*\otimes U^*\rangle_r:= \langle W,
W^*\rangle_r\langle U, U^*\rangle_r\,.
$$
As a consequence, the numbering of components in a tensor power
$V^{*\otimes k}$ is reverse to that in a tensor power $V^{\otimes
k}$
$$
V^{*\otimes k}:=V^*_k\otimes\dots V^*_2\otimes V^*_1\,, \qquad
V^{\otimes k}:=V_1\otimes V_2\otimes\dots \otimes V_k\,.
$$
This should be always kept in mind when working with operators
marked by numbers of spaces where these operators act (like in
formulae (\ref{YB}) and all other similar expressions).

Extend now the braiding (\ref{flip}) onto the space $V^*\otimes
V^*$. Below we show that, requiring a consistence of the extended
braiding and an invariance of the pairing (\ref{A1}), we have the
only choice \be R(x^i\otimes x^j) = x^r\otimes
x^s\,R^{ji}_{\;sr}\,.
\label{e-flip} \ee Therefore, by analogy with the construction of
Section \ref{sec:g-fo}, we can define the representations of the
Hecke algebras $H_k(q)$, $k\in {\Bbb K}$, in tensor powers
$V^{*\otimes k}$, construct the projectors $E^\lambda_a$ and
introduce subspaces $V^*_{(\lambda,a)} \subset V^{*\otimes k}$ as
images of the corresponding projectors (see
(\ref{Y-proj})--(\ref{v-lit-r})). Taking into account the above
remark on the numbering of tensor product components, one can
show
that for any Young tableau $(\lambda,a)$ there exists the only
tableau $(\lambda, a^\prime)$ such that the spaces
$V_{(\lambda,a)}$
and $V^*_{(\lambda,a^\prime)}$ are dual with respect to the form
(\ref{A1}).

By definition, the class of objects of the category $\swv$
consists
of all  direct sums of spaces $\vl\otimes V^*_\mu$ and
$V^*_\mu\otimes \vl$, where $\lambda$ and $\mu$ are partitions of
non-negative integers. The zero partition corresponds to the
basic
space $V_0:=V$ or to its dual space $V^*_0:=V^*$. The ground
field
${\Bbb K}$ is treated as the unit object of the category $\swv$
$$
{\Bbb K}\otimes V = V = V\otimes {\Bbb K}.
$$

Define now the class of morphisms of $\swv$. First of all, we
should
define the set of braiding morphisms $R_{U,W}$ realizing the
isomorphisms $U\otimes W \cong W\otimes U$ for any two objects
$U$
and $W$. The braidings $R_{V_\lambda,V_\mu}$ and
$R_{V^*_\lambda,V^*_\mu}$ are completely determined by $R_{V,V}$
and
$R_{V^*,V^*}$ given by (\ref{flip}) and (\ref{e-flip}).
Therefore,
we only need  consistent definitions of $R_{V,V^*}$ and
$R_{V^*,V}$, since the braidings $R_{V_\lambda,V^*_\mu}$ and
$R_{V^*_\lambda,V_\mu}$ can be then constructed by standard
methods
(cf, for example, \cite{GLS1}). The consistency condition is the
following requirement. Having defined the four braidings
mentioned
above, we get a linear operator on $(V\oplus V^*)^{\otimes 2}$.
Our
definitions are consistent if this operator satisfies the
Yang-Baxter equation. This problem is solved by the proposition
below, whose main idea belongs to V.Lyubashenko (cf. \cite{LS}
and
the references therein).

\begin{proposition}
\label{pr:6} Let $\Psi$ be the skew-inverse operator {\rm
(\ref{psi})} of the skew-invertible braiding $R$. Consider an
extension of $R$  to the linear operator
$$
R:(V\oplus V^*)^{\ot 2}\rightarrow (V\oplus V^*)^{\ot 2}
$$
(we keep the same notation for the extended operator) in
accordance
with the formulae {\rm \be
\begin{array}{lcl}
V\ot V^*\rightarrow V^*\ot V &: \quad &
R(x_i\ot x^j)=x^k\ot x_l\,(R^{-1})_{\;ki}^{lj}\,,\\
\rule{0pt}{6mm} V^*\ot V\rightarrow V\ot V^* &:\quad &
R(x^j\ot x_i)= x_k\ot x^l\,\Psi_{\;li}^{kj}\,,\\
\rule{0pt}{6mm} V^*\ot V^*\rightarrow V^*\ot V^* &:\quad &
R( x^i\ot x^j)=x^k\ot x^l\,R^{ji}_{\;lk}\,,\\
\rule{0pt}{6mm} V\ot V\rightarrow V\ot V &: \quad & R(x_i\ot
x_j)=x_k\ot x_l\,R_{\;ij}^{kl}\,.
\end{array}
\label{ext} \ee} Then the extended operator $R$ is a braiding,
i.e.
it satisfies the Yang-Baxter equation {\rm (\ref{YB})} on the
space
$(V\oplus V^*)^{\ot 3}$.
\end{proposition}

{\bf Proof.\ }  Since $R$ is a linear operator, it suffices to
prove the proposition on the basis vectors of the space $(V\oplus
V^*)^{\ot 3}$. This space splits into the direct sum of eight
subspaces (from $V\ot V\ot V$ to $V^*\ot V^*\ot V^*$) and the
verification of the statement of proposition on each of these
subspaces is a matter of straightforward calculations based on
formulae (\ref{ext}).\hfill \rule{6.5pt}{6.5pt}
\medskip

At the second step of our construction, we assume that a linear
combination, the product, the direct sum and the tensor product
of a
finite family of categorical morphisms is a morphism too.

Then, following \cite{T}, we require the morphisms to be
natural (or functorial). This means that
$$
(g\ot f)\circ R_{U,W}=R_{U',W'}\circ (f\ot g)
$$
where $f:U\rightarrow U'$ and $g:W\rightarrow W'$ are two
categorical morphisms. As a consequence, we get the necessary
condition for a map $f:U\rightarrow U'$ to be a categorical
morphism: \be ({\rm id}_W\ot f)\circ R_{U,W}=R_{U',W}\circ (f\ot
{\rm id_W}), \qquad (f\ot {\rm id}_W)\circ R_{W,U}=R_{W,U'}\circ
({\rm id}_W\ot f). \label{nec-con} \ee A map $f$ satisfying this
condition will be called {\em $R$-invariant}. Thus, any
categorical
morphism must be $R$-invariant.

\begin{proposition}
Provided $R$ satisfies {\rm (\ref{ext})}, the following claims
hold
true:
\begin{enumerate}
\item[{\rm 1.}]
The pairing {\rm (\ref{A1})} is $R$-invariant.

\item[{\rm 2.}]
The linear map $\pi_r:\K\to V^*\ot V$ generated by {\rm \be
1\stackrel{\pi_r}{\mapsto} \sum_{i=1}^Nx^i\ot x_i,\label{B1} \ee}
is
also $R$-invariant.
\end{enumerate}
\end{proposition}

\noindent {\bf Proof.\ } In proving  Claim 1, one can confine
oneself to considering the simplest case of the formula
(\ref{nec-con}) when $W=V$ or $W=V^*$. This is a consequence of
the
structure of objects of the category $\swv$. In other words, we
have
to show the commutativity of the diagram \be
\begin{array}{rcl}
(V\ot V^*)\ot V^\#&\;
\stackrel{(\ref{ext})}{\longleftrightarrow}\;&
V^\#\ot(V\ot V^*)\\
\rule{0pt}{7mm} \langle\,,\rangle_r\ot {\rm id} \;\downarrow
& & \downarrow\;{\rm id}\ot\langle\,,\rangle_r\\
\rule{0pt}{7mm} \K\ot V^\#& = &V^\#\ot \K
\end{array}
\label{r-inv} \ee where $V^\#$ stands for $V$ or $V^*$. The
commutativity of the diagram immediately follows from formulae
(\ref{ext}), definitions (\ref{A1}), (\ref{psi}) and the
definition
of the inverse matrix $R^{-1}$.

Next, the same reasoning shows, that Claim 2 is equivalent to
the commutativity of the diagram
$$
\begin{array}{rcl}
(V^*\ot V)\ot V^\#&\;
\stackrel{(\ref{ext})}{\longleftrightarrow}\;&
V^\#\ot(V^*\ot V)\\
\rule{0pt}{7mm} \pi_r\ot {\rm id} \;\uparrow
& & \uparrow\;{\rm id}\ot\pi_r\\
\rule{0pt}{7mm} \K\ot V^\#& = &V^\#\ot \K
\end{array}
$$
which can be proved similarly to the previous case.\hfill
\rule{6.5pt}{6.5pt}
\medskip
\begin{remark}
\rm Note, that the $R$-invariance of the maps (\ref{A1}) and
(\ref{B1}) is a motivation of the extension (\ref{ext}) of the
initial braiding $R$. It can be shown that such an extension is
unique.
\end{remark}

In what follows, besides the right form (\ref{A1}), we also
need a {\it left} nondegenerated bilinear form
$$
\langle\,,\rangle_l:\;V^*\ot V\to \K,
$$
with the additional requirement that the above pairing would be
$R$-invariant. This requirement prevents us from setting $\langle
x^i, x_j \rangle_l = \delta^i_j$, since it is not an
$R$-invariant
pairing (direct consequence of (\ref{ext})).

Let us choose the form $\langle\,,\rangle_l$ in such a way that
the
following diagram would be commutative \be
\begin{array}{rcl}
V^*\ot V &\;\stackrel{(\ref{ext})}{\longrightarrow}\;& V\ot V^*\\
\rule{0pt}{7mm} \langle\,,\rangle_l\;\downarrow
& & \downarrow\;\langle\,,\rangle_r\\
\rule{0pt}{7mm} \K & = & \K \,.
\end{array}
\label{l-form} \ee
A simple calculation based on
(\ref{l-form}), leads to the following explicit expression \be
\langle
x^i,x_j\rangle_l = B^i_{\;j}, \label{A2} \ee where the matrix
$\|B^i_{\;j}\|$ is defined in (\ref{def:B-C}). Such a choice
guarantees the $R$-invariance of the left pairing
$\langle\,,\rangle_l$. The commutativity of the corresponding
diagram (analogous to (\ref{r-inv})) can be easily verified
using of (\ref{ext}) and (\ref{dva}).
\begin{remark}
\rm Note, that the backward diagram \be
\begin{array}{rcl}
V^*\ot V &\;\stackrel{(\ref{ext})}{\longleftarrow}\;& V\ot V^*\\
\rule{0pt}{7mm} \langle\,,\rangle_l\;\downarrow
& & \downarrow\;\langle\,,\rangle_r\\
\rule{0pt}{7mm} \K & = & \K \,.
\end{array}
\label{l-form2} \ee is {\it not} commutative with the definition
(\ref{A2}). In a tensor category one can define the left pairing
in
such a way that both diagrams (\ref{l-form}) and (\ref{l-form2})
would be commutative, while in a quasitensor category it is
impossible. This is a consequence of the fact that the braiding
$R$
is not involutive: $R^2\not= I$.

In principle, we could demand the commutativity of the above
diagram
instead of diagram (\ref{l-form}). In this case in the right hand
side of (\ref{A2}) we would obtain an additional factor
$q^{2(m-n)}$. Actually, both variants are equivalent and choosing
between them is a matter of taste.
\end{remark}

Now we can find another basis $\{\lb{x}{i}\}_{1\le i\le N}$ of
$V^*$
which is dual to the basis $\{x_i\}_{1\le i\le n}$ with respect
to
the left form \be \lb{x}{i}:=q^{2(m-n)}\,C^i_{\;j}\,x^j, \quad
\Rightarrow\quad\langle\lb{x}{i},x_j\rangle_l = \delta^i_j.
\label{l-bas} \ee The normalizing factor in the definition of the
basis vector $\lb{x}{i}$ is chosen in accordance with Corollary
\ref{cor:4}.

So, we have two $R$-invariant bilinear forms and two basic sets
$\{x^i\}$ and $\{\lb{x}{i}\}$ in the space $V^*$ which are dual
to
the basis $\{x_i\}$ of the space $V$ with respect to the right
and
left forms correspondingly (see (\ref{A1}) and (\ref{l-bas})).
Due
to this reason, we refer to $\{x^i\}$ (resp.,
$\{\lb{x}{i}\}$) as the {\it right} (resp., {\it left})
basis of $V^*$.

Using of (\ref{dva}) one can rewrite the formulae (\ref{ext})
in terms of the left basis.
\begin{corollary}
In terms of the left basis $\{\lb{x}{i}\}_{1\le i\le N}$ of the
space $V^*$, the extension of the braiding $R$ defined by {\rm
(\ref{ext})} has the following form
$$
\begin{array}{l}
R(x_i\ot \lb{x}{j})=\lb{x}{k}\ot x_l\,\Psi_{\;ik}^{jl}\,,\\
\rule{0pt}{6mm} R(\lb{x}{j}\ot x_i)= x_k\ot \lb{x}{l}\,
(R^{-1})_{\;il}^{jk}\,,\\
\rule{0pt}{6mm} R(\lb{x}{i}\ot \lb{x}{j}) = \lb{x}{k}\ot
\lb{x}{l}
\,R^{ji}_{\;lk}\,,\\
\rule{0pt}{6mm}
R(x_i\ot x_j)=x_k\ot x_l\,R_{\;ij}^{kl}\,.\\
\end{array}
\eqno{(\ref{ext}\,')}
$$
Besides, the linear map $\pi_l:\K\to V\ot V^*$ generated by {\rm
\be
1\stackrel{\pi_l}{\mapsto} \sum_{i=1}^Nx_i\ot \lb{x}{i}
\label{B2}
\ee} is $R$-invariant.
\end{corollary}

{\bf Proof.} We prove the first formula in the above list
(\ref{ext}$\,'$), the others are proved in the same way. Taking
into account the definition of the left basis (\ref{l-bas}) and
the
first formula in (\ref{ext}), we get (recall the summation over
repeated indices)
$$
R(x_i\ot \lb{x}{j}) = x^u\ot x_l\, q^{2(m-n)}\,C^j_{\;s}
(R^{-1})^{ls}_{\;ui} = \lb{x}{k}\ot x_l\,q^{2(m-n)}
C^j_{\;s}(R^{-1})^{ls}_{\;ui} B^u_{\;k},
$$
where in the last equality we come back from the right basis to
the
left one by the formula inverse to (\ref{l-bas})
$$
x^u = B^u_{\;k}\,\lb{x}{k}.
$$
Then, from the second line formulae of (\ref{dva}) and Corollary
\ref{cor:4} we deduce
$$
q^{2(m-n)} C_{1}R^{-1}_{21}B_{2} = \Psi_{12}
$$
which allows us to make the following substitution in the above
line
of transformations
$$
q^{2(m-n)} C^j_{\;s}(R^{-1})^{ls}_{\;ui}B^u_{\;k} =
\Psi^{jl}_{\;ik}.
$$
So, we finally get
$$
R(x_i\ot \lb{x}{j}) = \lb{x}{k}\ot x_l\,\Psi^{jl}_{\;ik},
$$
that is the first line formula in (\ref{ext}$\,'$).

As for the $R$-invariance of the map $\pi_l$, it can be proved by
straightforward calculations on the base of (\ref{ext}) or
(\ref{ext}$\,'$) in the same way as it was done in proving of
(\ref{B1}). \hfill\rule{6.5pt}{6.5pt}
\medskip

Now we are able to define the categorical morphisms of $\swv$.
Together with the identical map, the list of the morphisms
includes
(\ref{A1}), (\ref{B1}), (\ref{A2}), (\ref{B2}) and all maps
(\ref{ext}) (or equivalently (\ref{ext}$\,'$)). Besides, as we
have
already mentioned above, any linear combination, the product
(successive application), the tensor product or the direct sum of
categorical morphisms is also a categorical morphism.

\begin{remark} {\rm Given a particular braiding $R$, one can in
principle compose a bigger list of $R$-invariant maps than that
mentioned above. Thus, for a super-space $V=V_0\oplus V_1$ the
projections $V\to V_0$ and $V\to V_1$ are $R$-invariant maps.}
\end{remark}

In what follows we are especially interested in the objects
$V^*\ot V$ and $V\ot V^*$ which are isomorphic to the space ${\rm
End}(V)$ of endomorphisms of the space $V$. Upon fixing the basis
$\{x_i\}$ in the space $V$, we come to the standard basis
$\{h_i^j\}
= x_i\otimes \lb{x}{j}$ in the space $V\otimes V^*$. Defining the
action of an element $v\ot v^*\in V\otimes V^*$ on a vector $u\in
V$
by the usual rule
$$
(v\ot v^*)( u):= v\langle v^*,u\rangle_l
$$
we get the action
$$
h_i^j(x_k) = \delta_k^j\,x_i
$$
and the multiplication table of elements $h_i^j$ treated as
endomorphisms of the space $V$
$$
h_i^j\circ h_k^s=\delta^j_k\,h_i^s.
$$

Fixing the right basis $\{x^i\}$ in $V^*$, we come to another
basis $l_i^j = x_i\otimes x^j$ of $V\otimes V^*$ with the
properties
(see (\ref{A2})) \be l_i^i(x_k) = B_{\;k}^j\,x_i,\qquad
l_i^j\circ
l_k^s = B_{\;k}^j\,l_i^s\,. \label{l-tu} \ee Taking into account
(\ref{l-bas}), we find the connection of the two basis sets \be
h_i^j = q^{2(m-n)}\,C^j_{\;k} l_i^k\,. \label{corr2} \ee

Introduce now the linear map $\Tr_R: \End(V)\to \Bbb K$ by means
of
the categorical morphism (\ref{A1}) \be \Tr_R(l^i_j)= \langle
x_j,x^i\rangle_r = \delta^i_j. \label{R-tr} \ee This map is
called the {\it $R$-trace} in what follows. By virtue of
(\ref{corr2}), the $R$-trace of an operator ${\sf F}\in \End(V)$
is
given by \be \Tr_R({\sf F}) = q^{2(m-n)}\,Tr(F\cdot C)\,,
\label{trF} \ee where $F$ is the matrix of the operator ${\sf F}$
with respect to basis $\{x_i\}$.

To complete the Section, we calculate the $R$-dimension of the
objects $\vl$ of our category. By definition, the $R$-dimension
of
an object $\vl\subset V^{\otimes k}$, $\lambda\vdash k$, is given
by
\be \dim_R\vl :=\Tr_R({\rm id}_{\vl}) = q^{2k(m-n)} Tr_{(1\dots
k)}(C_1\dots C_kE_a^\lambda)\,. \label{r-dim} \ee Basing on
(\ref{id-conn}), one can prove, that the above definition does
not
depend on the value of $a$. Besides, the $R$-dimension is an
ad\-di\-tive-mul\-ti\-pli\-ca\-tive functional as the classical
dimension
$$
\dim_R(U\otimes W) = \dim_RU\dim_RW\,,\qquad \dim_R(U\oplus W) =
\dim_RU + \dim_RW\,.
$$

Let us introduce the $R$-analogs $Q_\pm(t)$ of the HP series
$P_\pm(t)$ (\ref{hp-ser}) by the relation
$$
Q_\pm(t)=\sum_{k\ge 0} t^k\,\dim_R\Lambda^k_\pm(V)\,.
$$
Then the following proposition holds true.

\begin{proposition}
Given a skew-invertible Hecke symmetry with a bi-rank $(m|n)$, we
find the following properties of the series $Q_\pm$:
\begin{enumerate}
\item
if $m-n=0$, then $\dim_R\, \vl=0$ for any $\la\not=0$ and
therefore
$Q_+(t)=Q_-(t)= 1$.
\item If $m-n>0$, then
$$
\dim_R \vl = \dim_R \vl^* = s_\lambda (q^{m-n-1}, q^{m-n-3},...,
q^{1-m+n}),
$$
and therefore
$$
Q_-(t)=\sum_{k=0}^{m-n}{m-n \choose k }_q t^k\,, \qquad {p\choose
k}_q:=\frac{p_q(p-1)_q\dots (p-k+1)_q}{k_q(k-1)_q\dots 2_q1_q}\,.
$$

\item If $m-n<0$, then
$$
\dim_R \vl = \dim_R \vl^* = s_{\la^*}(q^{n-m-1}, q^{n-m-3},...,
q^{1-n+m}),
$$
where $\la^*$ is the conjugate partition, and therefore
$$
Q_+(t)=\sum_{k=0}^{n-m}{n-m \choose k }_q t^k\,.
$$
\end{enumerate}
\end{proposition}
\noindent{\bf Proof.} The proposition is proved by direct
calculations on the basis of definition (\ref{r-dim}). The
calculation are analogous to those in the even case (cf., for
example, \cite{GLS1}). \hfill\rule{6.5pt}{6.5pt}
\medskip

Emphasize, that $Q_\pm(t)$ depend only on the bi-rank of a given
Hecke symmetry $R$, whereas the corresponding HP series
$P_\pm(t)$,
besides the bi-rank, essentially depend on a concrete form of
$R$.

\section{mREA: definition and deformation properties}
\label{sec:mrea}

If $R$ is an involutive ($R^2=I$) skew-invertible symmetry, then
the
space $\End(V)$ can be endowed with the structure of {\em a
generalized Lie algebra} (cf. \cite{G1,G3}). The corresponding
enveloping algebra $U_R(\End(V))$ is defined as the following
quotient \be U_R(\End(V)) = T(\End(V))/\langle{\cal J}_R\rangle
\label{env} \ee where $\langle {\cal J}_R\rangle$ is a two sided
ideal of the free tensor algebra $T(\End(V))$ generated by the
subset ${\cal J}_R\subset T(\End(V))$ of the following form \be
{\cal J}_R = \{ X\ot Y-\Ren(X\ot Y)-X\circ Y+\circ \Ren(X\ot
Y)\;|\;\forall X,Y\in \End(V)\}\,. \label{g-Li} \ee Here $\circ$
is
the product in $\End(V)$ considered as an associative algebra of
linear operators on $V$ and the linear operator
$\Ren:\End(V)^{\ot
2}\to \End(V)^{\ot 2}$ is an extension of the braiding $R$ to the
space $\End(V)^{\ot 2}$. Its explicit form can be obtained
using (\ref{ext}).

Namely, choosing the basis $l^i_j=x_j\ot x^i$ in
the space $\End(V)$ and applying the corresponding formulae from
the
list (\ref{ext}) we find \be \Ren(l^i_j\ot
l^k_s)=l^{a_1}_{b_1}\ot
l^{a_2}_{b_2}\,
(R^{-1})^{b_2c_1}_{\;a_1c_2}\,R^{b_1c_2}_{\;jr_1}\,
R^{kr_2}_{\;a_2c_1}\,\Psi^{r_1i}_{\;r_2s}\,. \label{Rend} \ee In
order to present this formula in more transparent form we
introduce
the matrix notations which will be useful in what follows. Define
the $N\times N$ (recall, that $N=\dim V$) matrix $L$ with the
matrix
elements \be L_i^{\;j} = l^j_i, \label{l-mat} \ee where the first
(lower) index numerate rows and the second (upper) one numerates
columns of $L$. Then, introducing the matrix $\bar R$ transposed
to
$R$
$$
\bar R_{i_1i_2}^{\;\, j_1j_2} = R^{j_1j_2}_{\;\,i_1i_2}\,,
$$
we denote \be L_{\bar 1} = L\ot I, \quad L_{\bar 2} = \bar
R_{12}L_{\bar 1}\bar R^{-1}_{12}. \label{l-cop} \ee Now,
multiplying
the both sides of (\ref{Rend}) by $R$ and $R^{-1}$ and taking
into
account the definition of $\Psi$ (\ref{psi}), we represent
formula
(\ref{Rend}) in the equivalent form \be \Ren(L_{\bar 1}\ot
L_{\bar
2}) = L_{\bar 2}\ot L_{\bar 1}\,, \label{Rend2} \ee where the
summation over the corresponding matrix indices is understood.

Note, that a direct generalization of (\ref{env})--(\ref{g-Li})
with (\ref{Rend2}) from the involutive symmetry to the Hecke
symmetry case leads to an algebra which possesses bad deformation
properties and a poor representation theory. Fortunately, for any
skew-invertible Hecke symmetry $R$ there exists another
generalization of the enveloping algebra $U_R(\End(V))$
(\ref{env})
which has good deformation properties (see Proposition
\ref{pro:12} below) and coincides with the enveloping algebra
$U_R(\End(V))$  when $R$ is involutive.

\begin{definition}
\label{def:10} \rm The associative algebra generated by the unit
element $e_{\cal L}$ and the indeterminates $l_i^{j},\,\,1\le
i,j\le
N$ subject to the  system of relations \be R_{\;ij}^{kl}\,
l_k^m\,R_{\;ml}^{pq} \,l_p^r - l_i^a\,R_{\;aj}^{bc}\,
l_b^d\,R_{\;dc}^{rq} - \hbar(R_{\;ij}^{aq}\,l_a^r -
l_i^b\,R_{\;bj}^{rq}) = 0 \label{mREA1} \ee is called the
reflection equation algebra (REA) and denoted  $\lrq$ if $\h=0$,
and it is called the modified reflection equation algebra
(mREA) and denoted  $\lrqh$ if $\h\not=0$.
\end{definition}

The defining relations (\ref{mREA1}) can be presented in a
compact
form in terms of the matrix $L$ (\ref{l-mat}) and the transposed
matrix
$\bar R$ \be \bar R_{12} L_{1}\bar R_{12} L_{1} - L_{1}\bar
R_{12}
L_{1}\bar R_{12} - \hbar(\bar R_{12}\,L_{1} - L_{1}\,\bar R_{12})
=
0\,. \label{mREA} \ee
\begin{remark}
\label{rem:18} \rm Note, that by a linear transformation
of generators $l_i^j \mapsto m_i^j$ (at $q\not=\pm 1$)
$$
M = Ie_{\cal L} - \omega\hbar^{-1}L, \qquad \omega = q -
q^{-1}\,,
\quad M = \|m_i^j\|\,,
$$
we arrive at the following form of commutation relations
(\ref{mREA})
\be \bar R_{12} M_{1}\bar R_{12} M_{1} -
 M_{1}\bar R_{12} M_{1}\bar R_{12} = 0\,.
\label{Rmm} \ee This means, that the algebras ${\cal
L}(R_q,\hbar)$
and ${\cal L}(R_q)$ are isomorphic at $q\not=\pm 1$. The basis of
mREA generators with commutation relations (\ref{mREA}) is more
suitable for treating this algebra as an analog of universal
enveloping algebra $U(gl(m|n))$.
\end{remark}

Let us prove that the commutation relations (\ref{mREA}) are
consistent with the structure of the category $\swv$ in the
following sense. We treat the mREA as a quotient of the tensor
algebra $T(V\otimes V^*)$ over the two sided ideal, generated by
relations (\ref{mREA}) or, equivalently, by (\ref{Rmm}). These
relations are consistent with the structure of the category if
the
corresponding two sided ideal is invariant with respect to
braidings
of the category, or, in other words, if the mREA commutation
relations are {\it $R$-invariant}.

\begin{proposition}
\label{pro:eqvr} The commutation relations {\rm (\ref{mREA})} are
$R$-invariant.
\end{proposition}

\noindent{\bf Proof.} To prove the proposition, it is sufficient
to
show, that the commutation relations (\ref{mREA}) are preserved
when
commuting with $V$ or $V^*$ with respect to the braidings of the
category $\swv$. This can be done by straightforward calculations
on
the base of formulae (\ref{ext}) and $l_i^j = x_i\otimes x^j$. To
simplify the calculations, working with generators $m_i^j$
(\ref{Rmm}) is more convenient.

For example, taking a basis vector $x_i\in V$ we get, using
(\ref{ext})
$$
{\sf R}(x_{i_1}\otimes m_{i_2}^{j_2}) = \bar
R_{i_1i_2}^{\;a_1a_2}m_{a_1}^{b_1} (\bar
R^{-1})_{b_1a_2}^{\;c_1j_2}\otimes x_{c_1}\quad {\rm or}\quad
{\sf
R}(x_1\otimes M_2) = \bar R_{12}M_1\bar R^{-1}_{12}\otimes x_1\,,
$$
where ${\sf R}$ is a general notation for the corresponding
braiding, ${\sf R} = R_{V,V\otimes V^*}$ in the above formulae.
Now
we can directly get the desired result
\begin{eqnarray*}
x_1\otimes (\bar R_{23} M_{2}\bar R_{23} M_{2} - M_{2}\bar R_{23}
M_{2}\bar R_{23} ) &\stackrel{{\sf
R}}{\rightarrow} &\\
&&\hspace*{-10mm} \bar R_{12}\bar R_{23} (\bar R_{12} M_{1}\bar
R_{12} M_{1} -  M_{1}\bar R_{12} M_{1}\bar R_{12}) \bar
R^{-1}_{23}\bar R^{-1}_{12}\otimes x_1\,.
\end{eqnarray*}
The commutativity with $V^*$ is verified
analogously.\hfill\rule{6.5pt}{6.5pt}
\medskip

\begin{proposition} Let $R$ be an involutive skew-invertible
symmetry. Then the commutative relations among the generators
$\{l_i^j\}$ of the algebra $U_R(\End(V))$ {\rm (\ref{env})} are
equivalent to {\rm (\ref{mREA})} with $\h=1$. Therefore,
according
to Definition {\rm \ref{def:10}}, the algebra {\rm (\ref{env})}
coincides with mREA ${\cal L}_q(R,1)$.
\end{proposition}

\noindent{\bf Proof.} In involutive case $R=R^{-1}$ by
definition.
Therefore, the matrix $L_{\overline 2}$ defined in (\ref{l-cop})
can
be written as $L_{\overline 2} = \bar R_{12}L_1\bar R_{12}$. This
leads to the following action of $R_{\End}$ (\ref{Rend2})
$$
R_{\End}(L_1\otimes \bar R_{12}L_1\bar R_{12}) = \bar
R_{12}L_1\bar
R_{12}\otimes L_1.
$$
Now, by setting in (\ref{g-Li}) $X=L_1$, $Y=\bar R_{12}L_1\bar
R_{12}$, and taking into account the multiplication table for the
generators $l_i^j$ (\ref{l-tu}), we get
$$
X\circ Y = L_1R_{12}\,, \qquad \circ R_{\End}(X\otimes Y) =
R_{12}L_1\,,
$$
Together with the above form of the action of $R_{\End}$ this
allows us
to represent the set ${\cal J}_R$ (\ref{g-Li}) in the form
(\ref{mREA}) with $\hbar=1$. \hfill\rule{6.5pt}{6.5pt}
\medskip

The main deformation property of the mREA is given by the
following
proposition.

\begin{proposition}
\label{pro:12} Let $R$ be a skew-invertible involutive Hecke
symmetry: $R^2 = I$ and $U\subset {\Bbb K}$ be a neighbourhood of
$1\in {\Bbb K}$. Consider a family of skew-invertible Hecke
symmetries $R_q$, analytically depending on $q\in U$ and
satisfying
the condition $R_1 = R$. Denote the homogeneous component of
$\lrq$
of the $k$-th order by ${\cal L}^{(k)}(R_q)$. Then, provided $q$
is
generic, the following claims hold true.
\begin{enumerate}
\item[{\rm 1.}] \hspace*{10mm}
$ \dim\, {\cal L}^{(k)}(R_q)=\dim\, {\cal L}^{(k)}(R),\qquad
\forall\, k\ge 0\,. $
\item[{\rm 2.}]\hspace*{10mm}
$ {\rm Gr\,} \lrqh\cong \lrq\,, $\par where ${\rm Gr \,} \lrqh$
is
the graded algebra associated to the filtrated algebra $\lrqh$.
\end{enumerate}
\end{proposition}

\noindent {\bf Proof.} The verification of the item 1 is based on
the following observations. Below we construct a projector $({\rm
Span}(l_i^j))^{\otimes 3} \to {\cal L}^{(3)}(R_q)$. The explicit
form of the projector allows us to conclude that its rank is
constant for generic $q\in U$. Therefore \be \dim\,{\cal
L}^{(3)}(R_q)=\dim\, {\cal L}^{(3)}(R). \label{d3} \ee For an
involutive $R$, the algebra ${\cal L}(R)$ is the symmetric
algebra
of the linear space ${\rm Span}(l_i^j)$ equipped with the
involutive
braiding $R_{\End}$. This algebra is a Koszul one. (For the
definition of this notion the reader is referred to \cite{PP}.)
The
Koszul property of ${\cal L}(R)$ follows easily from exactness of
the Koszul complex of the second kind constructed in \cite{G3}.

Now we apply the result of \cite{PP} (generalizing \cite{Dr})
asserting
that the Koszul property of ${\cal L}(R)$ and relation (\ref{d3})
imply Claim 1 of our proposition. Moreover, it can be shown,
that for a generic $q\in U$ the algebra ${\cal L}(R_q)$ is also a
Kozsul algebra.

In order to prove Claim 2 of the proposition, we consider
the map $[\,\,,\,\,]$ sending the l.h.s. of (\ref{mREA1}) to its
r.h.s. As was shown in \cite{G4}, this map satisfies the Jacobi
relation in form of \cite{PP}. Then by  a generalization of
the PBW theorem, given in  \cite{PP} (cf., also, \cite{BG}), we
arrive at Claim 2.\hfill \rule{6.5pt}{6.5pt}
\medskip

\begin{remark}
\rm Note, that skew-invertible Hecke symmetries with
non-classical
HP series $P_-(t)$, constructed by methods of \cite{G3},
analytically depend on $q$ in a neighbourhood of 1.
\end{remark}

Now, we pass to a construction of the projector mentioned in the
proof
of proposition \ref{pro:12}. Rewrite REA ${\cal L}(R_q)$
((\ref{mREA}) at $\hbar =0$) in an equivalent form \be \bar
R_{12}L_{\bar 1}L_{\bar 2} - L_{\bar 1}L_{\bar 2} \bar R_{12} =
0.
\label{q-m-a} \ee Consider the unital associative algebra ${\frak
L}$ over $\Bbb K$ freely generated by $N^2$ generators $l_i^j$
$$
{\frak L} = {\Bbb K}\langle l_i^j\rangle\quad 1\le i,j\le N.
$$
The algebra ${\cal L}(R_q)$ is the quotient of $\frak L$ over the
two sided ideal $\langle {\cal I}_-\rangle$ generated by the left
hand side of (\ref{q-m-a}) \be {\cal L}(R_q) = \frak L/\langle
{\cal
I}_-\rangle, \qquad {\cal I}_- = L_{\overline 1}L_{\overline 2} -
\bar R_{12} L_{\overline 1}L_{\overline 2}\bar R_{12}^{-1}\,.
\label{REA-f} \ee

As a vector space, the algebra ${\frak L}$ can be decomposed into
a direct sum of homogeneous components
$$
{\frak L} = \bigoplus_{k\ge 0} {\frak L}_k, \quad {\frak
L}_0\cong
{\Bbb K},
$$
where each ${\frak L}_k$ is the linear span of the $k$-th order
monomials in generators $l_i^j$. The following basis turns out to
be
convenient to use \be {\frak L}_k = {\rm Span}[L_{\overline
1}L_{\overline 2}\dots L_{\overline k}]. \label{l-basis} \ee This
notation means that ${\cal L}_k$ is spanned by the matrix
elements
of the right hand side matrix. The matrices $L_{\overline m}$ are
defined by the recurrent rule \be L_{\bar 1} = L\otimes I, \qquad
L_{\overline{k+1}} = \bar R_{k}L_{\overline k}\bar R^{-1}_{k},
\quad
k\ge 1, \label{rec-r} \ee where the shorthand notation
$R_k:=R_{k\,k+1}$ will be systematically used below.
\begin{remark}
\label{rem-bas} \rm Note, that due to definitions (\ref{q-m-a}),
(\ref{rec-r}), and the Yang-Baxter equation for $R$ the following
relation holds \be \bar R_kL_{\overline k}L_{\overline{k+1}} =
L_{\overline k}L_{\overline{k+1}} \bar R_k\,\quad \forall\,k\ge
1.
\label{q-m-a-2} \ee Relation (\ref{q-m-a-2}) is typical for the
so-called quantum matrix algebras, REA being a particular case of
them. For detailed treatment of the question the reader is
referred
to \cite{IOP}.
\end{remark}

For the algebra ${\cal L}(R_q)$ we have an analogous vector space
decomposition
$$
{\cal L}(R_q) = \bigoplus_{k\ge 0}{\cal L}_k, \quad {\cal
L}_0\cong
{\Bbb K},\quad {\cal L}_k\subset {\frak L}_k.
$$
Let us try to describe the subspaces ${\cal L}_k$ explicitly. In
other words, we should find a series of projector operators
${\cal
S}_k: {\frak L}_k \rightarrow {\frak L}_k$ with the property
$$
{\rm Im}\,{\cal S}_k = {\cal L}_k\subset \frak L_3\,.
$$
Here we construct such projectors for the second and third order
components ${\cal L}_k$, $k=2,3$.

Introduce a linear operator $Q: {\frak L}_2\rightarrow {\frak
L}_2$
by the formula

\be Q(L_{\overline 1}L_{\overline 2}): = {\bar R}_1L_{\overline
1}L_{\overline 2}{\bar R}^{-1}_1 \ee or symbolically $Q = {\bar
R}_1\circ {\bar R}_1^{-1}$. Taking into account the Yang-Baxter
equation for $\bar R$ we can easily obtain that $Q$ also
satisfies
the Yang-Baxter equation \be Q_1Q_2Q_1 = Q_2Q_1Q_2\,,
\label{Q-yb}
\ee where $Q_1 = Q\otimes {\rm id}$ and $Q_2 = {\rm id}\otimes Q$
are obvious extensions of the operator $Q$ on the space $\frak
L_3$.
Moreover, using the fact that $R$ is a Hecke symmetry, one can
find
a minimal polynomial of the operator $Q$ \be (Q+q^2\,\mbox{\bf
I})(Q+q^{-2}\,\mbox{\bf I})(Q-\mbox{\bf I}) = 0, \quad \mbox{\bf
I}
:= I\circ I \label{Q-ch} \ee

As follows from (\ref{Q-ch}), the operator $Q$ has three
eigenvalues on the space $\frak L_2$. In obvious notations,
we get the vector space decomposition
$$
\frak L_2 = \frak L^{(-q^2)}\oplus \frak L^{(1)}\oplus \frak
L^{(-q^{-2})}\,.
$$
With the help of the Hecke condition (\ref{Hec}) one can find the
expressions for the corresponding projectors
\begin{eqnarray}
&&{\cal P}^{(-q^2)} = P_+(R)\circ P_-(R)\nonumber \\
\rule{0pt}{6mm}
&&{\cal P}^{(-q^{-2})} = P_-(R)\circ P_+(R)\nonumber\\
\rule{0pt}{6mm} &&{\cal P}^{(1)} =  P_+(R)\circ P_+(R)+
P_-(R)\circ
P_-(R)\,,\label{q-symtr}
\end{eqnarray}
where
$$
P_\pm(R) = {q^{\mp 1}I\pm\bar R\over 2_q}.
$$
Indeed, the direct calculation shows that
$$
Q{\cal P}^{(a)} = {\cal P}^{(a)}Q = a\,{\cal P}^{(a)}, \quad a=
-q^{\pm 2},\,1
$$
and, on the other hand, the operators ${\cal P}^{(a)}$ form the
complete
set of orthonormal projectors on $\frak L_2$
$$
{\cal P}^{(a)}{\cal P}^{(b)} = \delta^{ab}\,{\cal P}^{(a)}, \quad
{\cal P}^{(-q^2)} + {\cal P}^{(1)} + {\cal P}^{(-q^{-2})} =
\mbox{\bf I}.
$$

Now we take into account that (\ref{q-m-a}) is equivalent to
$$
(Q-\mbox{\bf I})(L_{\overline 1}L_{\overline 2}) = 0.
$$
This means that the second order component ${\cal L}_2$ of the
REA
${\cal L}(R_q)$ coincides (as a vector space) with the subspace
${\frak L}^{(1)}\subset {\frak L}_2$. Introduce a couple of the
orthonormal projection operators on  $\frak L_2$
$$
{\cal S}:= {\cal P}^{(1)},\quad {\cal A}:= {\cal P}^{(-q^2)}+
{\cal
P}^{(-q^{-2})}, \quad {\cal SA} = {\cal A\,S} = 0, \quad {\cal S}
+
{\cal A} = \mbox{\bf I}.
$$
These operators can be expressed in terms of $Q$
\begin{eqnarray}
&&{\cal S} = \frac{1}{2_q^2}\left((q^2+q^{-2})\, \mbox{\bf I}+Q +
Q^{-1}\right)\,,
\label{s-Q}\\
&&{\cal A} = \frac{1}{2_q^2}\left(2\,\mbox{\bf I}- Q -
Q^{-1}\right)\,,\nonumber
\end{eqnarray}
where the inverse operator $Q^{-1}$ can be obtained from
(\ref{Q-ch})
$$
Q^{-1} = Q^2+(q^2 -1 +q^{-2})Q-(q^2 -1 +q^{-2})\mbox{\bf I}\,.
$$
One can show that
$$
{\rm Span}({\cal I_-}) = {\rm Im}\,{\cal A}
$$
as vector subspaces in $\frak L_2$, where ${\cal I}_-$ is defined
in
(\ref{REA-f}).

The above considerations proves the following proposition.
\begin{proposition}
\label{prop:p2} The second order homogeneous component ${\cal
L}_2$
of the REA ${\cal L}(R_q)$ coincides with the image of the
$q$-symmetrizer ${\cal S}$ {\rm \be {\cal L}_2 = {\rm Im}\, {\cal
S}
= {\frak L}_2/{\rm Im}\,{\cal A}. \ee}
\end{proposition}

Let us pass to the third order homogeneous component $\frak L_3$
and
find the projector on the corresponding component ${\cal
L}_3\subset
\frak L_3$ of the REA ${\cal L}(R_q)$.

Extend the projectors $\cal S$ and ${\cal A}$ onto the subspace
$\frak L_3$. For the projector ${\cal S}$ the extension is given
by
the two operators ${\cal S}_1$ and ${\cal S}_2$ in accordance
with
the rule (see definition (\ref{q-symtr}))
$$
{\cal S}_1 :={\cal P}^{(1)}(R_1), \qquad {\cal S}_2 := {\cal
P}^{(1)}(R_2)\,,
$$
which means that
$$
{\cal S}_1(xyz) :=({\cal S}(xy))z, \qquad {\cal S}_2(xyz)
:=x({\cal
S}(yz)), \qquad \forall\, xyz\in \frak L_3\,.
$$
The formulae for $\cal A$  are analogous.

At this point we can see an advantage of using the basis
(\ref{l-basis}). Indeed, the {\it quadratic} homogeneous
component
${\cal L}_2$ can be embedded into $\frak L_3$ in different ways,
the
following two being the most important in the sequel
$$
{\cal L}_2\cdot {\frak L}_1\subset \frak L_3\quad {\rm and }\quad
{\frak L}_1\cdot{\cal L}_2\subset \frak L_3\,.
$$
As follows from (\ref{q-m-a}), (\ref{q-m-a-2}) and proposition
\ref{prop:p2}, these embeddings can be identified with the images
of
the operators ${\cal S}_1$ and ${\cal S}_2$
$$
{\cal L}_2\cdot{\frak L}_1 = {\cal S}_1(\frak L_3)\,, \qquad
{\frak
L}_1\cdot{\cal L}_2 = {\cal S}_2(\frak L_3).
$$
The following technical lemma plays a crucial role in the further
consideration.
\begin{lemma}
\label{lem:5} The $q$-symmetrizer ${\cal S}$ obeys the following
{\it fifth} order relation on the subspace $\frak L_3$: {\rm \be
{\cal S}_1{\cal S}_2{\cal S}_1{\cal S}_2{\cal S}_1 - a\,{\cal
S}_1{\cal S}_2{\cal S}_1 + b\,{\cal S}_1= {\cal S}_2{\cal
S}_1{\cal
S}_2{\cal S}_1{\cal S}_2 - a\,{\cal S}_2{\cal S}_1{\cal S}_2
+b\,{\cal S}_2 \label{5-rel} \ee} where {\rm \be a =
(q^4+q^2+4+q^{-2}+q^{-4})/2_q^4\quad b = 4_q^2/2_q^8\,.
\label{a-b}
\ee}
\end{lemma}

\noindent {\bf Proof.\ } The lemma is proved by a direct
calculation. The calculation can be considerably simplified if
for
${\cal S}$  one uses expression (\ref{s-Q}) instead of the
initial
definition (\ref{q-symtr}). \hfill \rule{6.5pt}{6.5pt}
\medskip

Consider the operator ${\cal S}^{(3)}:\frak L_3\to \frak L_3$
defined by \be {\cal S}^{(3)} = \frac{2_q^6}{4\cdot
3_q^2}\,\left({\cal S}_1{\cal S}_2{\cal S}_1{\cal S}_2{\cal S}_1
-a\,{\cal S}_1{\cal S}_2{\cal S}_1 +b\,{\cal S}_1\right)\,,
\label{3-symtr} \ee with $a$ and $b$ given in (\ref{a-b}). Due to
(\ref{5-rel}) there exists an equivalent form of the above
operator
\be {\cal S}^{(3)} = \frac{2_q^6}{4\cdot 3_q^2}\,\left({\cal
S}_2{\cal S}_1{\cal S}_2{\cal S}_1{\cal S}_2 -a\,{\cal S}_2{\cal
S}_1{\cal S}_2 +b\,{\cal S}_2\right)\,. \label{3-symtr2} \ee In
fact, the operator ${\cal S}^{(3)}$ is the projector on ${\cal
L}_3\subset \frak L_3$ we are looking for.

\begin{proposition}
The third order homogeneous component ${\cal L}_3$ of the REA
${\cal
L}(R_q)$ is the image of the projection operator ${\cal S}^{(3)}$
under its action on $\frak L_3$
$$
{\cal L}_3 = {\rm Im}\,{\cal S}^{(3)}\,.
$$
\end{proposition}

\noindent {\bf Proof.\ } The fact that $({\cal S}^{(3)})^2 =
{\cal
S}^{(3)}$ can be verified by a straightforward calculation.

Consider now the projection of the relation (\ref{REA-f}) onto
the
third order homogeneous component: \be {\cal L}_3 = \frak
L_3/\langle{\cal I_-}\rangle_3, \qquad \langle{\cal I_-}\rangle_3
=
\frak L_1\cdot {\rm Im}\,{\cal A}_2 \cup {\rm Im}\,{\cal A}_1
\cdot
\frak L_1\,. \label{3-hcomp} \ee As can be seen from
(\ref{3-symtr})
and (\ref{3-symtr2}), $\langle{\cal I_-}\rangle_3 \subseteq {\rm
Ker}\, {\cal S}^{(3)}$ and, therefore, ${\rm Im}\,{\cal
S}^{(3)}\subseteq {\cal L}_3$.

>From the other side, since ${\cal S} + {\cal A} = \mbox{\bf I}$,
the
subspace ${\cal L}_3$ given in (\ref{3-hcomp}) can be presented
as
$$
{\cal L}_3 = \frak L_1\cdot {\rm Im}\, {\cal S}_2\cap {\rm
Im}\,{\cal S}_1\cdot \frak L_1\,.
$$
Comparing this form of ${\cal L}_3$ with the structure of ${\cal
S}^{(3)}$ given in (\ref{3-symtr}) and (\ref{3-symtr2}), we get
${\rm Im}\,{\cal S}^{(3)}\supseteq {\cal L}_3$. This relation
together with the opposite inclusion obtained above completes the
proof. \hfill \rule{6.5pt}{6.5pt}

\section{The braided bialgebra structure and representation
theory} \label{sec:7}

In this Section we  consider finite dimensional representations
of the mREA (\ref{mREA}) in the category $\swv$. Note, that the
class of finite dimensional representations of mREA is wider, for
instance, it includes a large number of one-dimensional
representations. For the particular case of the $U_q(sl(m))$
$R$-matrix all such a representations were classified in
\cite{Mu1}.
Besides, the finite dimensional representations of mREA can be
constructed on the base of the $U_q(sl(m))$ representations,
since
the REA ${\cal L}(R_q)$ can be embedded (as an algebra) into the
quantum group in this particular case.

The representation theory developed below does not depend on a
particular choice of the $R$-matrix and works well in the general
situation, when the quantum group does not exist. Besides, an
important property of the suggested theory is the equivariance of
the representations we are dealing with. By definition, a
representation $\rho_U$ of mREA in a space $U$ is called {\it
equivariant} if the map
$$
\End(V) \rightarrow \End(U):\quad l_i^j\mapsto \rho_U(l_i^j)
$$
is a morphism of the category $\swv$.

This property has an important consequence, which will be used
below. Namely, given an mREA module $W$ with equivariant
representation $\rho_W: {\cal L}(R_q,1)\rightarrow \End(W)$, the
diagram \be
\begin{array}{ccc}
U\otimes (A\otimes W)&\stackrel{{\sf R}}{\longleftrightarrow}
&(A\otimes W)\otimes U\\
\rule{0pt}{6mm} \hspace*{10mm} \downarrow{\rm id}\otimes \rho_W &
&\hspace*{-10mm}
\rho_W\otimes {\rm id}\downarrow\\
\rule{0pt}{6mm} U\otimes W&\stackrel{{\sf
R}}{\longleftrightarrow} &
W\otimes U
\end{array}
\label{eqvr} \ee is commutative for any object $U$ of the
category
$\swv$ and for any subspace $A\subset {\cal L}(R_q,1)$. The
equivariance condition allows us to define an mREA representation
in
the tensor product of mREA modules.

For the particular case of an even Hecke symmetry of the rank
$(m|0)$, the equivariant representation theory of the associated
mREA turns out to be similar to the representation theory of the
algebra $U(sl(m))$. At the beginning of the Section \ref{sec:SWH}
we
mentioned the specific peculiarity of the even case. Namely, in
the
corresponding category of the mREA representations, the space
$V^*$
can be identified with the object $V_{(1^{m-1})}$ and for
constructing the complete representation theory it suffices, in
fact, to define the mREA-module structure on the space $V$ and on
its tensor powers $V^{\otimes k}$. Any tensor product $V^{\otimes
k}$ is a reducible mREA-module and expands into the direct sum
(\ref{V-decom}) of mREA-invariant subspaces $V_\lambda$ (for
details
cf. \cite{GS2,S}).

However, the construction of the cited papers is insufficient for
the
treatment of the general case of the bi-rank $(m|n)$. The reason
is
that in general case we have to construct representations in the
tensor products  $V^{*\otimes k}$ independently of those in
tensor
products $V^{\otimes k}$ and the central problem here consists in
extending the mREA-module structure on the tensor product of
modules
$V_\lambda\otimes V_\mu^*$.

In the present Section we suggest the regular procedure for
constructing the mREA representations  which works well
independently of the bi-rank of the Hecke symmetry. Our
construction
is based on the braided bialgebra structure in the mREA.
Throughout
this Section we set $\hbar = 1$ in the mREA commutation relations
(\ref{mREA}).

The main component of the braided bialgebra is the coproduct
$\Delta$, which is a homomorphism of the mREA ${\cal L}(R_q,1)$
into
an associative braided algebra $\mbox{\bf L$(R_q)$}$ which is
defined as follows.

\begin{itemize}
\item As a vector space over the field $\Bbb K$ the algebra
$\mbox{\bf L$(R_q)$}$ is isomorphic to the tensor product of two
copies of mREA
$$
\mbox{\bf L$(R_q)$} = {\cal L}(R_q,1) \otimes {\cal L}(R_q,1)\,.
$$
\item The product $\star : \mbox{\bf L$(R_q)$}^{\otimes
2}\rightarrow \mbox{\bf L$(R_q)$}$ is defined by the rule \be
(a_1\otimes b_1)\star (a_2\otimes b_2):=a_1 a^\prime_2 \otimes
b^\prime_1 b_2\,,\qquad a_i\otimes b_i \in \mbox{\bf L$(R_q)$}\,,
\label{br-pr} \ee where $a_1a^\prime_2$ and $b_1b^\prime_2$ are
the
usual product of elements of mREA, while $a^\prime_1$ and
$b^\prime_1$ result from the action of the braiding $R_{\rm End}$
(see
(\ref{Rend2})) on the tensor product $b_1\otimes a_2$ \be
a^\prime_2\otimes b^\prime_1:= R_{\rm End}(b_1\otimes a_2)\,.
\label{ti-el} \ee
\end{itemize}

We should verify that product (\ref{br-pr}) is indeed
associative.
For this we need the following lemma.
\begin{lemma}
\label{lem:kp} Consider copies of the matrix $L$ defined in {\rm
(\ref{rec-r})}. Then the following relation holds {\rm \be
\Ren(L_{\overline k}\otimes L_{\overline p}) = L_{\overline
p}\otimes L_{\overline k}\,, \qquad \forall\,k< p,\; k,p\in {\Bbb
N}\,. \ee }
\end{lemma}

\noindent {\bf Proof.} The proof consists in a straightforward
calculation on the base of relation (\ref{Rend2}) rewritten in a
slightly modified form
$$
\Ren(L_1{\bar R}_{12}\otimes L_1) = {\bar R}_{12}L_1{\bar
R}^{-1}_{12}\otimes L_1 {\bar R}_{12}\,,
$$
and on the Yang-Baxter equation (\ref{YB}) which allows one to
interchange the chains of $R$-matrices forming the copies
$L_{\overline k}$ and $L_{\overline p}$.\hfill\rule{6.5pt}{6.5pt}
\smallskip

Now the associativity of (\ref{br-pr}) can be easily proved for
the
elements $X^{i}_{r,s}\in \mbox{\bf L$(R_q)$}$ whose components
are
homogeneous monomials in generators of ${\cal L}(R_q,1)$:
$$
X^{i}_{r,s}:=L_{\overline{i}}\dots L_{\overline{i+r-1}}\otimes
L_{\overline{i+r}}\dots L_{\overline{i+r+s-1}}\,,
$$
where we represent the homogeneous components of mREA as a linear
span of elements similar to those (\ref{l-basis}).

Note, that for any triple $X_{r_1,s_1}^{i_1}$,
$X_{r_2,s_2}^{i_2}$
and $X_{r_3,s_3}^{i_3}$ we can always take
$$
i_3\ge i_2+r_2+s_2 \ge i_1+r_1+s_1
$$
and then the associativity condition
$$
(X_{r_1,s_1}^{i_1}\star X_{r_2,s_2}^{i_2})\star X_{r_3,s_3}^{i_3}
=
X_{r_1,s_1}^{i_1}\star (X_{r_2,s_2}^{i_2}\star X_{r_3,s_3}^{i_3})
$$
is an immediate consequence of lemma \ref{lem:kp}. Since any
element
of $\mbox{\bf L$(R_q)$}$ can be presented as a linear combination
of
some $X^{(i)}_{r,s}$ we conclude that (\ref{br-pr}) defines an
associative product in $\mbox{\bf L$(R_q)$}$.

Note, that the mREA are isomorphic to two subalgebras of
$\mbox{\bf
L$(R_q)$}$ by  the following embeddings
$$
a\mapsto e_{\cal L}\otimes a\quad {\rm or}\quad a\mapsto a\otimes
e_{\cal L}\,,
$$
where $e_{\cal L}$ is the unit of mREA ${\cal L}(R_q,1)$. This
can
be easily obtained from the fact that the unit element $e_{\cal
L}$
trivially commutes with any $a\in {\cal L}(R_q,1)$ with respect
to
the braiding $\Ren$. As a consequence we have
$$
(e_{\cal L}\otimes a_1)\star (e_{\cal L}\otimes a_1) = (e_{\cal
L}\otimes a_1a_2)\quad {\rm and} \quad (a_1\otimes e_{\cal
L})\star
(a_2\otimes e_{\cal L}) = (a_1a_2\otimes e_{\cal L})\,.
$$

Define a linear map $\Delta: {\cal L}(R_q,1)\rightarrow \mbox{\bf
L$(R_q)$}$ by the following rules: \be
\begin{array}{l}
\Delta(e_{\cal L}):= e_{\cal L}\otimes e_{\cal L}\\
\rule{0pt}{6mm} \Delta(l_i^j) := l_i^j\otimes e_{\cal L}+e_{\cal
L}\otimes l_i^j -
(q-q^{-1})\sum_k l_i^k\otimes l_k^j\\
\rule{0pt}{6mm}
\Delta(ab):=\Delta(a)\star\Delta(b)\qquad\forall\,a,b\in {\cal
L}(R_q,1)\,.
\end{array}
\label{copr} \ee In addition to (\ref{copr}), we introduce a
linear
map $\varepsilon :{\cal L}(R_q,1) \rightarrow {\Bbb K}$ \be
\begin{array}{l}
\varepsilon(e_{\cal L}):= 1\\
\rule{0pt}{6mm}
\varepsilon(l_i^j):= 0\\
\rule{0pt}{6mm} \varepsilon(ab):=
\varepsilon(a)\varepsilon(b)\qquad
\forall\,a,b\in {\cal L}(R_q,1)\,.
\end{array}
\label{coed} \ee

The proposition below establishes an important property of the
maps
$\Delta$ and $\varepsilon$.

\begin{proposition}
\label{pro:28} The maps $\Delta$ and $\varepsilon$ given in {\rm
(\ref{copr})} and {\rm (\ref{coed})} are respectively the
coproduct
and counit of the braided bialgebra structure on the mREA ${\cal
L}(R_q,1)$.
\end{proposition}

\noindent {\bf Proof.} First, we prove that the map $\Delta$
defines
an algebra homomorphism ${\cal L}(R_q,1)\rightarrow \mbox{\bf
L$(R_q)$}$. It is convenient to work with the generators $m_i^j$
introduced in Remark \ref{rem:18}. In terms of these generators
the
map $\Delta$ reads \be \Delta(m_i^j) = \sum_s m_i^s\otimes
m_s^j\,,
\label{co-mm} \ee or, in the matrix form
$$
\Delta(M_{\overline k}) = M_{\overline k}\otimes M_{\overline k}
\qquad \forall\, k\ge 1\,.
$$
Taking  definitions (\ref{br-pr}) and (\ref{Rend2}) into account
we
find
$$
\Delta(M_{\overline 1}M_{\overline 2}) = \Delta(M_{\overline
1})\star \Delta(M_{\overline 2}) = M_{\overline 1}M_{\overline 2}
\otimes M_{\overline 1}M_{\overline 2}\,.
$$
Comparing this with (\ref{Rmm}) we finally get
$$
R_{12}\Delta(M_{\overline 1}M_{\overline 2})  =
\Delta(M_{\overline
1}M_{\overline 2}) R_{12}
$$
which means that the map $\Delta$ is an algebra homomorphism.
Note,
that braided coproduct (\ref{co-mm}) was suggested in \cite{M}.

The interrelation between $\varepsilon$ and $\Delta$
$$
({\rm id}\otimes \varepsilon)\,\Delta = {\rm id} =
(\varepsilon\otimes {\rm id})\, \Delta
$$
is verified trivially. \hfill\rule{6.5pt}{6.5pt}
\medskip

Consider now the question: which representations of $\mbox{\bf
L$(R_q)$}$ can be constructed using mREA
representations? Given two equivariant mREA modules $U$ and $W$
with
representations $\rho_U:{\cal L}(R_q,1)\rightarrow \End(U)$ and
$\rho_W:{\cal L}(R_q,1)\rightarrow \End(W)$, we construct the map
$\rho_{U\otimes W}:\mbox{\bf L$(R_q)$}\rightarrow \End(U\otimes
W)$
by the following rule \be \rho_{U\otimes W}(a\otimes
b)\triangleright (u\otimes w) = (\rho_U(a)\triangleright
u^\prime)\otimes (\rho_W(b^\prime)\triangleright w)\,, \qquad
a\otimes b\in \mbox{\bf L$(R_q)$}\,, \label{rep-ll} \ee where the
symbol $\triangleright$ stands for the action of the
corresponding
operator, the element(s) $b^\prime$ and vector(s) $u^\prime$
result
from the action of the corresponding braiding (depending on $b$
and
$u$) of the category $\swv$
$$
u^\prime \otimes b^\prime:={\sf R}(b\otimes u)\,.
$$
By Proposition \ref{pro:eqvr} the definition
(\ref{rep-ll}) is self-consistent since the map $b\mapsto
\rho_W(b^\prime)$ is also a representation of mREA ${\cal
L}(R_q,1)$.

\begin{proposition}
\label{pro:29} The action {\rm (\ref{rep-ll})} defines a
representation of the algebra $\mbox{\bf L$(R_q)$}$.
\end{proposition}

\noindent{\bf Proof.} Consider two arbitrary elements of the
algebra
$\mbox{\bf L$(R_q)$}$
$$
X_i = (a_i\otimes b_i)\in\mbox{\bf L$(R_q)$},\quad i=1,2\,.
$$
We have to prove that $\forall\,u\otimes w\in U\otimes W$ the
following relation holds \be \rho_{U\otimes W}(X_1\star
X_2)\triangleright (u\otimes w) = \rho_{U\otimes
W}(X_1)\triangleright(\rho_{U\otimes W}(X_2)\triangleright
(u\otimes
w))\,. \label{doc-rep} \ee The left and right hand sides of
(\ref{doc-rep}) are actually the maps sending the element
$$
(a_1\otimes b_1)\otimes (a_2\otimes b_2)\otimes (u\otimes w)
$$
into a vector of the space $U\otimes W$. We prove that the
results of applying these maps to the above element are the same.

Let us introduce the shorthand notations
$$
{\sf R}(b_1\otimes a_2) = a^\prime_2\otimes b^\prime_1\,,\quad
{\sf
R}(b_2\otimes u) = u^\prime\otimes b^\prime_2\,,\quad {\sf
R}(b^\prime_1\otimes u^\prime) = u^{\prime\prime} \otimes
b^{\prime\prime}_1\,.
$$
The definitions (\ref{br-pr}) and (\ref{rep-ll}) allow us to
represent the left hand side of (\ref{doc-rep}) as the
composition
of the following morphisms
\begin{eqnarray*}
(a_1&\hspace*{-2.5mm}\otimes& \hspace*{-2.5mm} b_1)\otimes
(a_2\otimes b_2)\otimes (u\otimes w) \mapsto (a_1\otimes
a^\prime_2)\otimes (b^\prime_1\otimes b_2)\otimes (u\otimes w)
\mapsto (a_1a^\prime_2\otimes b^\prime_1b_2)\otimes (u\otimes
w)\\
&&\mapsto (a_1a^\prime_2\otimes
u^{\prime\prime})\otimes(b^{\prime\prime}_1b^\prime_2\otimes w)
\mapsto (\rho_U(a_1a^\prime_2)\triangleright
u^{\prime\prime})\otimes
(\rho_W(b^{\prime\prime}_1b^\prime_2)\triangleright w)\,.
\end{eqnarray*}
Now, take into account the equivariance condition (\ref{eqvr})
for
the representations $\rho_U$ and $\rho_W$. This condition means
that
under the action of the categorical braidings a vector
$\rho_U(a)\triangleright u$ commutes with any object in the same
manner as the element $a\otimes u$ does. Therefore the right hand
side of (\ref{doc-rep}) can be represented as the composition of
the
maps
\begin{eqnarray*}
(a_1&\hspace*{-2.5mm}\otimes& \hspace*{-2.5mm} b_1)\otimes
(a_2\otimes b_2)\otimes (u\otimes w) \mapsto (a_1\otimes
b_1)\otimes
(a_2 \otimes u^\prime)\otimes
(b^\prime_2\otimes w)\\
&&\mapsto (a_1\otimes b_1)\otimes (\rho_U(a_2)\triangleright
u^\prime\otimes \rho_W(b^\prime_2)\triangleright w)\mapsto
(a_1\otimes \rho_U(a^\prime_2)\triangleright
u^{\prime\prime})\otimes
(b^{\prime\prime}_1\otimes \rho_W(b^\prime_2)\triangleright w)\\
&&\mapsto
(\rho_U(a_1)\triangleright\rho_U(a^\prime_2)\triangleright
u^{\prime\prime})\otimes
(\rho_W(b^{\prime\prime}_1)\triangleright
\rho_W(b^\prime_2)\triangleright w) \mapsto
(\rho_U(a_1a^\prime_2)\triangleright u^{\prime\prime})\otimes
(\rho_W(b^{\prime\prime}_1b^\prime_2)\triangleright w)
\end{eqnarray*}
So, having started from the same initial element, the maps in the
left and right hand sides of (\ref{doc-rep}) give the same
resulting
vector in $U\otimes W$. Therefore, these maps are identical.
\hfill\rule{6.5pt}{6.5pt}
\medskip

\begin{corollary}
Let $U$ and $W$ be two ${\cal L}(R_q,1)$-modules with equivariant
representations $\rho_U$ and $\rho_W$. Then equivariant
representation ${\cal L}(R_q,1)\rightarrow \End(U\otimes W)$ is
given by the rule {\rm \be a\mapsto \rho_{U\otimes
W}(\Delta(a))\,,
\qquad \forall\,a\in {\cal L}(R_q,1)\,, \label{rep-prod} \ee }
where
the coproduct $\Delta$ and the map $\rho_{U\otimes W}$ are given
respectively by formulae {\rm (\ref{copr})} and {\rm
(\ref{rep-ll})}.
\end{corollary}
\noindent{\bf Proof.} This corollary is a direct consequence of
Propositions \ref{pro:28} and \ref{pro:29}.
\hfill\rule{6.5pt}{6.5pt}
\medskip

As was mentioned at the beginning of this Section, the
equivariant
representations of mREA in spaces $V_\lambda$, $\lambda\vdash
k\in
{\Bbb N}$, were constructed in \cite{GS2,S}. By the same method
we
can define the mREA representations in spaces $V^*_\lambda$.
Then,
the formula (\ref{rep-prod}) allows us to define the mREA-module
structure on any object of the category $\swv$.

To complete the picture, we shortly outline the main ideas of
\cite{GS2,S}, and, besides, prove the equivariance of
representation
in the space $V^*$. Contrary to the even case, for the Hecke
symmetry of the general bi-rank the equivariance of
representations
in the dual spaces $V^*_\lambda$ should be established
independently
of that in the spaces $V_\lambda$.

The basic left representation of mREA ${\cal L}(R_q,1)$ in the
space
$V$ is defined in terms of matrix elements of the operator $B$
\be
\rho_1(l_i^j)\triangleright x_k=B^j_{\;k}x_i. \label{l-rep} \ee
As
was shown in \cite{S}, the map $\rho_2:\lrqo\to \End(V^{\otimes
2})$
defined by
$$
\rho_2(l_i^j)\triangleright (x_{k_1}\ot x_{k_2})= (\rho_1(l_i^j)
\triangleright x_{k_1})\ot
x_{k_2}+\Bigl(R^{-1}\circ(\rho_1(l_i^j)\ot I) \circ R^{-1}\Bigr)
\triangleright (x_{k_1}\ot x_{k_2})
$$
is a representation. An extension of the basic representation up
to
the higher representations $\rho_p:{\cal L}(R_q,1)\to \End
(V^{\otimes p})$, $p\ge 3$, is defined in a similar way. It can
be
shown by direct calculations, that these extensions coincide with
the universal recipe (\ref{rep-prod}).

The representations of the above form are completely reducible
--- the space $V^{\otimes p}$ expands into the direct sum of
invariant subspaces $\vl$ labelled by partitions $\lambda\vdash
p$.
The restriction of the representation $\rho_p$ to subspaces $\vl$
is
obtained by the action of orthogonal projectors $E^\lambda_a$
(see
(\ref{Y-proj}) and (\ref{V-decom})) \be \rho_{\lambda,a} =
E^\lambda_a\circ\rho_p\circ E^\lambda_a\,, \label{restr} \ee the
modules with different $a$ being equivalent.

The basic representation $\rho_1^*:{\cal L}(R_q,1)\rightarrow
\End(V^*)$ is given by \be \rho_1^*(l_i^j)\triangleright x^k =
-x^r\,\bar R_{ri}^{\;kj}\,. \label{duals} \ee To prove the
equivariance of this representation, we need the following lemma.

\begin{lemma} Let $R$ be a skew-invertible Hecke symmetry.
Then the map {\rm \be V\ot V^*\to V^* \ot V\,:\quad x_i\ot
x^j\mapsto x^k\otimes x_l R_{\;ki}^{lj} \label{trans} \ee} is a
categorical morphism.
\end{lemma}
{\bf Proof.\ } We use the fact that for a Hecke symmetry $R =
R^{-1}
+(q-q^{-1})I$ (see (\ref{Hec})). Substituting this into
(\ref{trans})
$$
x_i\otimes x^j\mapsto x^k\otimes x_l\, (R^{-1})^{lj}_{\;ki}
+(q-q^{-1})\,\delta_i^j\,x^k\otimes x_k
$$
we find that the map in question is a linear combination of a
categorical morphism from the list (\ref{ext}) and of the map
$$
x_i\otimes x^j\mapsto \delta_i^j\,x^k\otimes x_k.
$$
It can be presented as a composition of categorical morphisms
$$
x_i\otimes x^j\stackrel{\langle\,,\,\rangle_r}{\mapsto}
\delta_i^j\,1\stackrel{\pi_r}{\mapsto} \delta_i^j\,x^k\otimes
x_k,
$$
where the categorical morphisms (\ref{A1}) and (\ref{B1}) were
used
in the consecutive order.

So, we conclude that the initial map (\ref{trans}) is a linear
combination of categorical morphism, therefore it is a
categorical
morphisms by definition. \hfill\rule{6.5pt}{6.5pt}
\medskip

\begin{proposition}
The representation {\rm (\ref{duals})} of the algebra $\lrqo$ in
the
space $V^*$ is equivariant.
\end{proposition}
{\bf Proof.\ } To prove the equivariance of $\rho_1^*$, we are to
show that the map $\rho_1^*: \lrqo\to \End(V^*)$ is a categorical
morphism.

Identifying $l_i^j$ with $x_i\otimes x^j$, we can treat any left
equivariant action of $l_i^j$ on the basis vector $x^k\in V^*$ as
a
categorical morphism
$$
V\otimes V^*\otimes V^* \to V^*.
$$
Construct such an action as the following composition of
morphisms
$$
V\otimes V^*\otimes V^* \stackrel{(\ref{trans})\otimes
I}{\longrightarrow} V^*\otimes V\otimes V^* \stackrel{I\otimes
\langle\,,\,\rangle_r}{\longrightarrow} V^*\otimes {\Bbb K}\cong
V^*
$$
which gives explicitly
$$
l_i^j\triangleright x^k = R^{kj}_{\;ri}\,x^r = x^r\bar
R^{\;kj}_{ri}.
$$
Up to a sign, this categorical morphism coincides with the left
representation (\ref{duals}).\hfill \rule{6.5pt}{6.5pt}
\medskip

Let us consider a particular example of the "adjoint" mREA
representation acting in the linear span of generators $l_i^j$.
Since
$$
{\rm Span}(l_i^j)\cong V\otimes V^*\,,
$$
the representation involved is constructed by the general formula
(\ref{rep-prod}) which now reads
$$
l_i^j\mapsto \rho_{V\otimes V^*}(\Delta(l_i^j))\,,
$$
where  we should take
the basic representations (\ref{l-rep}) and (\ref{duals})
as $\rho_{V}(l_i^j)$ and $\rho_{V^*}(l_i^j)$,
respectively. Omitting straightforward calculations, we write the
final result in the compact matrix form \be \rho_{V\otimes
V^*}(L_{\overline 1})\triangleright L_{\overline 2} = L_1R_{12} -
R_{12}L_1. \label{ad} \ee Applying the coproduct $\Delta$
(\ref{copr}), we  extend this representation onto any homogeneous
component of the mREA.

Note, that above action (\ref{ad}) is an $L$-linear part of the
defining commutation relations of the mREA (\ref{mREA}), if we
rewrite them in the equivalent form
$$
L_{\overline 1}L_{\overline 2} - R_{12}^{-1} L_{\overline
1}L_{\overline 2} R_{12} = L_1R_{12} - R_{12}L_1\,.
$$
In this sense the action (\ref{ad}) is similar to the adjoint
action
of a Lie algebra $\frak g$ on its universal enveloping algebra
$U(\frak g)$, which is also determined by the linear part of the
Lie
bracket and then is extended from the Lie algebra onto higher
components of $U(\frak g)$ by means of the standard coproduct
operation.

To end the Section, consider the question of the
"$sl$-reduction",
that is, the passing from mREA ${\cal L}(R_q,1)$ to the quotient
algebra \be {\cal SL}(R_q):= {\cal L}(R_q,1)/\langle
\Tr_RL\rangle\,, \qquad \Tr_RL:= Tr(CL). \label{sl-quo} \ee The
elements $\ell:=\Tr_RL$ is central in mREA, which can be easily
proved by calculating the $R$-trace in the second space from the
matrix relation (\ref{mREA}). In so doing, the formulae
(\ref{trRC})
is useful.

To describe the quotient algebra ${\cal SL}(R_q)$ explicitly, we
pass to the new set of  generators $\{f_i^j,\ell\}$, connected
with
the initial one by the linear transformation: \be l_i^j = f_i^j +
(Tr(C) )^{-1}\delta_i^j\,\ell\quad {\rm or}\quad L=F+(Tr(C)
)^{-1}I\,\ell\,, \label{shift} \ee where $F=\|f_i^j\|$.
Obviously,
$\Tr_RF = 0.$ Note, that by (\ref{trBC}), the above shift
is possible iff $p\not=r$.

In terms of new generators, the commutation relations of mREA
read
$$
\left\{
\begin{array}{l}
\displaystyle \bar R_{12}F_1\bar R_{12}F_1 -  F_1\bar
R_{12}F_1\bar
R_{12} = (e_{\cal L} - \frac{\omega}{Tr(C)}\,\ell)(\bar R_{12}F_1
-
F_1\bar R_{12})\\
\rule{0pt}{6mm} \ell \,F = F\,\ell\,,\qquad \Tr_RF = 0\,,
\end{array}
\right.
$$
where $\omega = q-q^{-1}$. Now, the quotient (\ref{sl-quo}) is
easy
to describe. The matrix $F=\|f_i^j\|$ of ${\cal SL}(R_q)$
generators
satisfy the same commutation relations (\ref{mREA}) as the matrix
$L$ \be \bar R_{12}F_1\bar R_{12}F_1 -  F_1\bar R_{12}F_1\bar
R_{12}
= \bar R_{12}F_1 - F_1\bar R_{12}\,, \quad \Tr_RF = 0\,,
\label{sl-rea} \ee but the generators $f_i^j$ are linearly
dependent
due to the relation $\Tr_RF = Tr(CF)=0$.

It is not difficult to rewrite the representation (\ref{ad}) in
terms of generators $f_i^j$ and $\ell$. Taking  relation
(\ref{shift}) into account, we find, after a short calculation
\begin{eqnarray}
&&\rho_{V\otimes V^*}(\ell)\triangleright \ell = 0,\qquad
\rho_{V\otimes V^*}(F_1)\triangleright \ell = 0\,,\nonumber\\
&&\rho_{V\otimes V^*}(\ell)\triangleright F_1 =
-\omega\,Tr(C)\,F_1\nonumber\\
&&\rho_{V\otimes V^*}(F_{\overline 1})\triangleright F_{\overline
2}
= F_1\bar R_{12} - \bar R_{12}F_1 + \omega \bar R_{12}F_1\bar
R_{12}^{-1}\,. \label{sl-ad}
\end{eqnarray}
Note, that relation (\ref{sl-ad}) defines the "adjoint"
representation of the quotient algebra ${\cal SL}(R_q)$, but,
contrary to the mREA ${\cal L}(R_q,1)$, this representation is
not
given by the linear part of the quadratic-linear commutation
relations
(\ref{sl-rea}).

Generally, given a representation $\rho:{\cal
L}(R_q,1)\rightarrow
\End(U)$ such that the element $\ell$ is a multiple of the unit
operator (for example, an irreducible representation)
$$
\rho(\ell) = \chi \,I_U\,, \qquad \chi\in{\Bbb K}\,
$$
we can construct the corresponding representation
$\tilde\rho:{\cal
SL}(R_q) \rightarrow \End(U)$ by the formula \cite{S} \be
\tilde\rho(f_i^j) =  \frac{1}{\xi}\left(\rho(l_i^j) -
(Tr(C))^{-1}\rho(\ell)\,\delta_i^j \right)\,, \qquad \xi =
1-(q-q^{-1})(Tr(C))^{-1}\chi\,. \label{sl-red} \ee At last, we
note,
that REA (\ref{Rmm}) admits a series of automorhisms $M\mapsto
zM$
with nonzero $z\in{\Bbb K}$. At the level of mREA representations
these automorphisms read (recall, that  $\hbar = 1$)
$$
\rho_U(l_i^j)\mapsto \rho_U^{z}(l_i^j) = z\rho_U(l_i^j) +
\delta_i^j(1-z)(q-q^{-1})^{-1}\,I_U\,.
$$
Using (\ref{sl-red}), one can show that the corresponding
representation $\tilde\rho_U$ of the algebra ${\cal SL}(R_q)$
constructed from $\rho_U^z$ does not depend on $z$, in other
words,
the whole class of mREA representations $\rho_U^z$ connected by
the
above automorphism gives the same representation of the quotient
algebra ${\cal SL}(R_q)$.

\begin{remark}
\rm In this connection we would like to discuss the problem of a
suitable definition of braided (quantum, generalized) Lie
algebras.
For the first time such an object was introduced in \cite{G1} as
a
data $(\gg,\,\sigma,\,[\,,\,])$ where $\gg$ is a vector space,
$\sigma:\gg^{\ot2}\to \gg^{\ot 2}$ is an involutive symmetry, and
$[\,,\,] :\gg^{\ot2}\to \gg$ is an operator ("braided Lie
bracket")
such that
\begin{enumerate}
\item
$[\,,\,]\sigma=-[\,,\,]$
\item
$\sigma [\,,\,]_{23}= [\,,\,]_{12}\sigma_{23} \sigma_{12}$
\item
$[\,,\,][\,,\,]_{23}(I+\sigma_{12}\sigma_{23}+\sigma_{23}\sigma_{
12})=0.$
\end{enumerate}
Remark that the third relation can be presented as follows \be
[\,,\,] [\,,\,]_{12}=[\,,\,][\,,\,]_{23}(I-\sigma_{12})
\label{genLie} \ee

A typical example is
$$
\gg=\End(V)\,,\quad \sigma=R_{\End}\,,\quad
[\,,\,]=\circ(I-\sigma)
$$
(in the setting of  Section 5). Another example can be obtained
by restricting the above operators onto the subspace of traceless
elements of the algebra $\End(V)$. The enveloping algebras of the
both braided Lie algebras can be defined by (\ref{env}).

Now, observe that relation (\ref{genLie}) takes the form
(\ref{ad})
if we put \be \gg={\rm Span}(l_i^j)\,,\quad \sigma(L_{\overline
1}
L_{\overline 2})=R_{12}^{-1} L_{\overline 1} L_{\overline 2}
R_{12}\,,\quad [L_{\overline 1},L_{\overline 2}]= L_{\overline
1}R_{12}- R_{12} L_{\overline 1}\,. \label{g-Lie} \ee

So, if we define a braided Lie algebra with such $\gg$, $\sigma$
and
$[\,,\,]$, the third axiom of the above list (in the form
(\ref{genLie}) is satisfied. By contrast, the relations from
the items 1 and 2 fail and must be modified. Thus, an analog of
the
item 1  can be presented in the form
$$
[\,,\,] {\cal S} =0\,,
$$
where $\cal S$ is defined in (\ref{s-Q}). The verification of
this
relation is straightforward and is left to the reader. In the
item
2, the map $\sigma$ must be replaced by $R_{\End}$. This is
consequence of the fact that the bracket $[\,,\,]$ in
(\ref{g-Lie})
is a categorical morphism. But if we restrict ourselves to the
traceless part of the space $\gg$, the relation (\ref{genLie})
fails
too.

So, it is somewhat contradictory to define a braided Lie algebra
in
the space $\End(V)$ (where the space $V$ is equipped with a
skew-invertible Hecke symmetry) with the use of the above three
axioms. However, in many papers (cf. \cite{Wo,GM}) the braided
(quantum) Lie algebras related to non-involutive braidings are
introduced via just these axioms or their slight modifications.
Taking in  consideration our observation (see Introduction) on
"braided Lie algebras"  related to braiding of the
Birman-Murakami-Wenzl type, we can conclude that there is no
"braided Lie algebras" satisfying the above list of axioms and,
at
the same time, such that their enveloping algebras
possess good  deformation property.
\end{remark}

\section{Quantization with a deformed trace}
\label{sec:tr-qua}

In this Section we consider the semiclassical structures arising
from the  mREA  $\lrqh$ (\ref{mREA}), provided that $R$ is the
standard $U_q(sl(m))$ Hecke symmetry (\ref{Rmat}). In this case
the
mREA $\lrqh$ is treated as a two-parameter deformation of the
commutative algebra $\K[gl(m)^*]$. We clarify the role of the
corresponding Poisson brackets in defining the quantum
homogeneous
spaces. At the end of the Section we study the infinitesimal
counterpart of the deformed $R$-trace.

Given the Hecke symmetry $R$ (\ref{Rmat}), we can find the
Poisson
pencil which is the semiclassical counterpart of the
two-parameter
algebra $\lrqh$ by a straightforward calculation. Indeed, setting
$q=1$ in (\ref{Rmat}), we pass from $\lrqh$ to the algebra
$U(gl(m)_\h)$. Therefore, the Poisson bracket corresponding to
the
deformation described by the parameter $\h$ is the linear
Poisson-Lie bracket on the space of functions on $gl(m)^*$.

In order to find the second generating bracket of the Poisson
pencil, we put $\h=0$ in (\ref{mREA}) coming thereby to the
non-modified REA. Introducing the matrix ${\cal R} = RP$
$$
{\cal R}= \sum_{i,j}^mq^{\delta_{ij}} h_i^i\otimes h_j^j +
(q-q^{-1})\sum_{i<j}^m h_i^j\otimes h_j^i\,,
$$
we transform the commutation relations of the REA into the form
\be
\bar{\cal R}_{12}\,L_{1}\,\bar{\cal R}_{21}\,L_{2} -
L_{2}\,\bar{\cal R}_{12}\,L_{1}\, \bar{\cal R}_{21} = 0\,,
\label{REA} \ee we recall that the bar over the symbol of a
matrix
means transposition.

On setting $q=e^{\nu}$, $\nu\in{\Bbb K}$, and noting, that ${\cal
R}=I$ at $q=1$, we come to the following expansion of ${\cal R}$
into a series in $\nu$: ${\cal R}=1+\nu\,\textsf{r}+ O(\nu^2)$,
where \be \textsf{r} = \sum_{i=1}^mh_i^i\otimes h_i^i +
2\sum_{i<j}^mh_i^j\otimes h_j^i. \label{s-class} \ee is the
classical $sl(m)$ $\textsf{r}$-matrix
$$
[\textsf{r}_{12},\textsf{r}_{13}] +
[\textsf{r}_{12},\textsf{r}_{23}]
+[\textsf{r}_{13},\textsf{r}_{23}]
= 0.
$$
Now, the part of the commutation relation (\ref{REA}) which is
linear in $\nu$ represent the second Poisson bracket in ${\Bbb
K}[gl(m)^*]$ \be \{L_1,L_2\}_r = L_{2}L_{1}\bar \textsf{r}_{21} -
\bar \textsf{r}_{12} L_{1}L_{2} + L_{2}\bar \textsf{r}_{12}L_{1}
-
L_{1}\bar \textsf{r}_{21}L_{2}\,. \label{pb-quadr} \ee This
formula
is defined on the matrix elements $l_i^j$ of the matrix $L$ which
form a basis of linear functions on the space $gl(m)^*$. The
extension of bracket (\ref{pb-quadr}) from the generators $l_i^j$
to
arbitrary functions (polynomials in generators) is described in
terms of vector fields on $gl(m)^*$. In order to obtain such an
extension, we introduce the matrices
$$
\textsf{r}_{\pm} = \frac{1}{2}\,(\textsf{r}_{12} \pm
\textsf{r}_{21}).
$$
As directly follows from (\ref{s-class}), the above matrices are
images of \be r_- = \sum_{i<j}^m(e_i^j\otimes e_j^i -
e_j^i\otimes
e_i^j):=\sum_{i<j}^me_i^j\wedge e_j^i \,,\qquad r_+=\sum_{i,j}^m
e_i^j\otimes e_j^i\,,\qquad r_\pm\in gl(m)^{\otimes 2}
\label{r-class} \ee under the fundamental vector representation
$e_i^j\mapsto h_i^j$, the elements $e_i^j$ being the standard
basis
of $gl(m)$
$$
[e_i^j,e_k^r] = \delta_k^j\,e_i^r - \delta_i^r\,e_k^j.
$$

Consider now the actions of $gl(m)$ on ${\Bbb K}[gl(m)^*]$ by the
left, right and adjoint vector fields \be e_i^j\triangleright
l_k^s:=\delta_k^j\,l_i^s\,,\qquad l_k^s\triangleleft
e_i^j:=\delta_i^s\,l_k^j\,,\qquad {\rm ad}\,e_i^j(l_k^s) =
e_i^j\triangleright l_k^s - l_k^s\triangleleft e_i^j\,,
\label{vec-ac} \ee which are extended on any polynomial from
${\Bbb
K}[gl(m)^*]$ by the Leibnitz rule. Then, taking into account the
above definitions of $r_\pm$ we can rewrite (\ref{pb-quadr}) in
the
general form \be \{f,g\}_r= \circ\, r_+^{\rm l,r}(f\ot g) -
\circ\,
r_+^{\rm r,l}(f\ot g) - \circ \,r_-^{\rm ad,ad}(f\ot g),\qquad
\forall\, f,g\in {\Bbb K}[gl(m)^*]. \label{brr} \ee Here $\circ:
{\Bbb K}[gl(m)^*]^{\otimes 2}\to {\Bbb K}[gl(m)^*]$ stands for
the
commutative pointwise product of functions on $gl(m)^*$ and
superscripts of $r_\pm$ denote the following actions \be r_-^{\rm
ad,ad}(f\ot g):=\sum_{i<j}^m{\rm ad}\,{e_i^j}(f) \wedge {\rm ad}
\,e_j^i(g), \quad r_+^{\rm l,r}(f\ot
g):=\sum_{i,j}(e_i^j\triangleright f) \otimes (g \triangleleft
e_j^i)\,. \label{r-sup} \ee Note that the brackets
$$
\{f,\,g\}_-=\circ\, r_-^{\rm ad,ad}(f\ot g)   \quad {\rm
and}\quad
\{f,g\}_+= \circ \,r_+^{\rm l,r}(f\ot g)-\circ\, r_+^{\rm
r,l}(f\ot
g)
$$
are {\em not} Poisson, since they do not obey the Jacobi
identity.

The bracket $\{\,\}_+$ is $gl(m)$-covariant, that is
$$
\ad X(\{f,g\}_+)=\{\ad X(f),g\}_++\{f,\ad X(g)\}_+,\qquad \forall
f,g\in \K[gl(m)^*],\quad X\in gl(m).
$$

The bracket (\ref{brr}) restricts onto any $gl(m)$-orbit ${\cal
O}\subset gl(m)^*$, since for any $f\in I_{\cal O}$, where
$I_{\cal
O}$ is an ideal of functions vanishing on this orbit, and for any
$g\in \K[gl(m)^*]$ we have
$$
\{f,g\}_r\in I_{\cal O}\,.
$$
This property is evident for the component $\{\,\}_-$, since it
is
defined via the adjoint vector fields. The proof for the
component
$\{\,\}_+$ is given in \cite{D}. In particular, the bracket
(\ref{brr}) can be restricted on the variety $c_1:=\sum l_i^i =
c$
where $c\in \K$ is a constant. On setting $c=0$, we get a Poisson
bracket on the algebra  $\K[sl(m)^*]$.

\begin{remark}
\rm Being restricted on $\K[\gg^*]$ where $\gg = sl(m)$, the
bracket
$\{\,\}_+$ admits the following interpretation \cite{G4,D}.
Consider
the space $\gg^{\ot 2}$ as an adjoint  $\gg$-module. It
decomposes
into the direct sum of submodules
$$
\gg^{\otimes 2} = \gg_s\oplus \gg_a,
$$
where $\gg_{s}(\gg_{a})$ is the symmetric (skew-symmetric)
subspace
of $\gg^{\ot 2}$. For $n>2$ there exist subspaces $\gg_+\subset
\gg_s$ and $\gg_-\subset \gg_a$ which are isomorphic to $\gg$
itself
as adjoint $\gg$-modules. Therefore there is a unique (up to a
factor) nontrivial $\gg$-covariant morphism $\beta:\gg^{\ot
2}\to\gg^{\ot 2}$ sending $\gg_-$ to $\gg_+$.

Since the space of linear functions on $\gg^*$ with the
Poisson-Lie
bracket is isomorphic to $\gg$ as a Lie algebra, the morphism
$\beta$ can also be defined on the whole algebra $\K[\gg^*]$ (via
the Leibnitz rule) and it gives a bracket which coincides (up to
a
factor) with $\{\,\}_+$. Note that for $n=2$ the component
$\gg_+$
vanishes and therefore the bracket $\{\,\}_+$ vanishes as well.
\end{remark}

\begin{proposition} The bracket {\rm (\ref{brr})} is compatible
with the Poisson-Lie one (their Schouten bracket vanishes) and,
therefore, any bracket of the pencil {\rm (\ref{pp})} is Poisson.
\end{proposition}
This claim is an immediate corollary of the fact that the family
of
algebras $\lrqh$ is a two parametric deformation of the
commutative
algebra ${\Bbb K}[gl(m)^*]$, but it can be also verified by
direct
calculations \cite{D}.

Now, let us consider the case $n=2$ in more detail. As we have
noticed
above, in this case the component $\{\,\}_{+}$ of the Poisson
bracket $\{\,\}_{r}$ vanishes and we have
$\{\,\}_{r}=\{\,\}_{-}$.
Let $\{H,E,F\}$ be the Cartan-Chevalley generators of $sl(2)$.
Then
$r_-= E\wedge F$ (see (\ref{r-class})) and the Poisson bracket
(\ref{brr}) reads
$$
\{a,b\}_r= - \ad E(a)\,\ad F(b) + \ad F(a)\,\ad E(b).
$$
Take the generators $\{e,f,h\}$ of $\K[sl(2)^*]$ which correspond
to the Cartan-Chevalley generators under the isomorphism of
$sl(2)$
and the Lie algebra of linear functions on $sl(2)^*$. A simple
calculation on the base of (\ref{vec-ac}) gives
$$
\{h,e\}_r= - 2eh\,,\quad \{h,f\}_r= 2fh\,,\quad \{e,f\}_r= -
h^2\,.
$$
Note, that this differs from the Poisson-Lie bracket only by the
factor $-h$ and, therefore,  each leaf of the bracket $\{\,\}_r$
lies in a leaf of the Poisson-Lie bracket. Moreover, the element
$c_2=h^2/2 +2ef$ is central with respect to the both brackets and
hence the corresponding Poisson pencil can be restricted onto the
quotient $\K[sl(2)^*]/\langle c_2-c\rangle$ for any $c\in \K$.

Let $\K=\C$ and elements
$$
x={1\over 2}\,(e-f),\quad y= {i\over 2}\,(e+f),\quad z={i\over
2}\,h
$$
be the generators of $su(2)$. In these generators we get the
Poisson
pencil $\{\,\}_{PL,r}$ where
$$
\begin{array}{lclcl}
\{x,y\}_{PL}=z,&\quad&\{y,z\}_{PL}=x,&\quad&\{z,x\}_{PL}=y,\\
\rule{0pt}{5mm} \{x,y\}_{r}=z^2,&\quad &\{y,z\}_{r}=xz,&\quad
&\{z,x\}_{r}=yz.
\end{array}
$$
Here we have renormalized the bracket $\{\,\}_r$ since this does
not
affect the Poisson pencil. The quadratic central element takes
the
form $c_2 = x^2+y^2+z^2$.

A particular bracket of this Poisson pencil (namely,
$\{\,\}_{PL}-\{\,\}_{r}$) appeared in \cite{Sh} (see Appendix by
J.-H.Lu and A.Weinstein) in studying a semiclassical counterpart
of
the quantum sphere.  In the cited paper, the quantum sphere was
represented as an operator algebra. This approach is based on the
work \cite{P}, where the quantum sphere was treated to be a
$C^*$-algebra and its irreducible representations (as
$C^*$-algebra)
were classified. Finally, in \cite{Sh} the quantum sphere was
presented in terms of functional analysis.

Our method of constructing the quantum sphere (or quantum
hyperboloid what is the same over the field $\K=\C$) is
completely
different. First of all, we quantize the KKS
bracket\footnote{Recall, that KKS bracket is a restriction of the
Poisson-Lie bracket on a coadjoint orbit of the Lie group in the
space dual to its Lie algebra.} on the sphere and realize the
resulting quantum algebra as the quotient \be U(su(2)_\h)/\langle
x^2+y^2+z^2 -c\rangle\,,\quad c\in {\Bbb K}\,,\;c\not =0.
\label{U-quot} \ee We are interested in finite dimensional
representations of this algebra. There exists a set of negative
values $c^{(k)}=-\h^2k(k+2)/4$, $k\in {\Bbb N}$, of the parameter
$c$ such that the quotient algebra (\ref{U-quot}) admits a finite
dimensional representation if $c=c^{(k)}$ for some $k\in {\Bbb
N}$.

Returning to the generators of the algebra $sl(2)$, we obtain a
one-parameter family of algebras
$$
\slch= U(sl(2)_\h)/\langle {1\over 2}\,H^2+EF+FE-c\rangle\,.
$$
Any algebra $\slch$ of this family as well as
$\slc=\K[sl(2)^*]/\langle{h^2\over 2}+ef+fe-c\rangle$, being
equipped with the $sl(2)$-action, can be expanded into a
multiplicity free direct sum of $sl(2)$-modules
$$
\slc \cong \bigoplus_{k\ge 0} V_k\,.
$$
Let $ \alpha : \slc\to \slch$ be an $sl(2)$-invariant map sending
the highest weight elements $e^{\ot k}\in \slc$ to $E^{\ot k}\in
\slch$. This requirement defines the map $\alpha$ completely.
Now,
the commutative algebra $\slc$ can be equipped with a new
noncommutative product $\star$ coming from the algebra $\slch$
\be
f\star_\h g=\al^{-1}(\al(f)\circ \al(g)),\qquad f,g\in\slc
\label{star-pr} \ee where $\circ$ is the product in the algebra
$\slch$. Thus, we have the quantized KKS bracket on the
hyperboloid
$c_2=c$ in the spirit of deformation quantization scheme ---
introducing a new product in the initial space of commutative
functions. Note, that according to \cite{R}, our algebraic
quantization cannot be extended on the function space
$C^\infty[S^2]$.

Now, deform  the algebra $\K[sl(2)^*]$ in $q$ and $\h$
"directions"
simultaneously. We get the mREA (\ref{mREA}) with $R$-matrix
given
in (\ref{Rmat}) where we should set $m=2$. Extracting the
$R$-traceless elements from the set of four mREA generators, we
come
to a unital associative algebra, generated by three linearly
independent elements $\{\hat H,\hat E,\hat F\}$ subject to the
system of commutation relations
$$
\begin{array}{l}
q^2 \hat H \hat E - \hat E \hat H = 2_q\h \, \hat E\,,\\
\rule{0pt}{7mm}
\hat H \hat F-q^2 \hat F \hat H = -2_q\h\,\hat F\,,\\
\rule{0pt}{7mm} \displaystyle q(\hat E \hat F-\hat F \hat E) =
\hat
H\Bigl(\h - \frac{(q^2-1)}{2_q}\, \hat H\Bigr)\,.
\end{array}
$$
Denote this algebra by $\slqh$. The element $C_q={\hat H^2\over
2_q}+q^{-1} \hat E \hat F+q\hat F\hat E$ is central and is called
{\em the braided Casimir}. Let us put
$$
\slcqh= \slqh/\langle C_q-c \rangle.
$$
We call this algebra the quantum hyperboloid or (considering it
over
the field $\K=\C$) the quantum sphere. It is a two-parameter
deformation of the initial commutative algebra $\slc$. For a
generic
value of $q$ it is possible to define a map $\alpha_q: {\cal
SL}^c\rightarrow \slcqh$ similar to the map $\alpha$ (but without
the equivariance property) and represent the product in $\slc$ in
the spirit of relation (\ref{star-pr}).

As in the case of algebra (\ref{U-quot}), there exists a series
of
values $c=c_k$, such that the corresponding quotient algebra
${\cal
SL}^{c_k}(q,\hbar)$ has a finite dimentional equivariant
representation. Its construction is described in Section
\ref{sec:7}. When $q\rightarrow 1$, we get a representation of
the
algebra ${\cal SL}^{c_k}(\hbar)$. By contrast, the representation
theory of the quantum sphere suggested  in \cite{P} has nothing
in
common with the theory of finite dimensional representations of
$sl(2)$ (or $su(2)$).

In general, by quantizing the KKS bracket on a semisimple orbit
we
represent the quantum algebra as an appropriate quotient of the
enveloping algebra $U(\gh)$ with $\gg=gl(m)$ or $sl(m)$. (Note,
that
if such an orbit is not generic the problem of finding defining
relations of the corresponding "quantum orbit" is somewhat
subtle,
cf. \cite{DM}). Finally, we additionally deform this quotient in
"q-direction"  and get some quotient of the algebra $\lqh$.

Observe, that on a generic orbit in $\gg^*$ where $\gg$ is a
simple
Lie algebra there exists a  family of nonequivalent Poisson
brackets, giving rise to $U_q(\gg)$-covariant algebras. One of
them
is the reduced Sklyanin bracket. It is often described in terms
of
the Bruhat decomposition (cf. \cite{LW}). The classification of
all
these brackets and their deformation quantization  are given in
\cite{DGS} and \cite{D}. The reduced Sklyanin bracket can be also
quantized in terms of the so-called Hopf-Galois extension (cf.
\cite{DGH}). But only the bracket (\ref{brr}), restricted to a
semisimple orbit, is compatible with the KKS bracket and the
quantization of the corresponding Poisson pencil  can be realized
in
the spirit of affine algebraic geometry, i.e.  via generators and
relations between them.

Note, that on the sphere (hyperboloid) the reduced Sklyanin
bracket
coincides with one of the bracket from the Poisson pencil
$\brkksr$
(it is also true for any symmetric orbit). So, it can be
quantized
via different approaches. However, for $m>2$ and for higher
dimensional orbits the notion of quantum orbits should be
concretized. It essentially depends on the bracket to be
quantized.

As for the other classical simple algebras $\gg$ of the $B$, $C$
or
$D$ series, there is no two-parameter deformation of the algebra
$\K[\gg^*]$ (cf. \cite{D}). Though a quadratic-linear algebra
(similar to $\lrqh$) can be constructed in this case (cf.
\cite{DGG}
for detail), we note that neither this algebra nor the associated
quadratic algebra is a deformation of its classical counterpart.

Let us complete this Section by considering a semiclassical
analog
of the quantum trace in the spirit of \cite{G2}. In that paper
 Poisson pencils similar to the above ones were considered but
 they were
generated by triangular classical $r$-matrices (they give rise to
involutive braidings). The main difference is that the result  of
the "double quantization" of Poisson pencil from \cite{G2} was
treated as the enveloping algebra of a generalized Lie algebra
and
its finite dimensional representations formed a tensor (not a
quasitensor) category.

As is known, on any symplectic variety there is a Liouville (or
invariant, or symplectic) measure $d\mu$ with basic property
$\int\{f,g\}d\mu=0$. In the framework of the deformation
quantization this measure gives rise to a trace with usual
properties (cf. \cite{GR}). It is just the case of the KKS
bracket
on a semisimple orbit. For non-symplectic Poisson brackets one
usually tries to describe its symplectic leaves and to quantize
them
separately, i.e., to associate  an operator algebra to each of
the
leaves. In the framework of our approach we are not dealing with
quantizing leaves of the bracket $\brr$ or  any bracket from the
Poisson pencil $\brkksr$ but we quantize this Poisson pencil as a
whole. In other words, we simultaneously q-deform all algebras
arising from "$\h$-quantization"  and arrive to operator algebras
with deformed traces.

Consider the Poisson pencil $\brkksr$ on a semisimple orbit
${\cal
O}\subset su(m)^*$. The bracket $\{\,,\}_r$ is not symplectic,
therefore the pencil involved has no Liouville measure on the
whole
orbit (a similar case was considered in \cite{G2}). Nevertheless,
the following proposition holds true independently of the
concrete
form of the matrix $r$.

\begin{proposition}
Let $\brkksr$ be the Poisson pencil on a semisimple orbit ${\cal
O}\subset\gg^*$, where $\gg=su(m)$ (or its complexification) and
$d\mu$ is the Liouville measure for the bracket $\brkks$. Then
the
quantity {\rm \be \langle a,b \rangle=\int_{\cal O}
\{a,b\}_r\,d\mu
\label{cocyc} \ee} is a cocycle with respect to the bracket
$\brkks$, i.e.
$$
\langle a,\{b,c\}_{KKS}\rangle+\langle
b,\{c,a\}_{KKS}\rangle+\langle c,\{a,b\}_{KKS}\rangle=0.
$$
\end{proposition}

This statement is a simple consequence of the fact that the
brackets
$\{\,,\}_{KKS}$ and $\{\,,\}_r$ are compatible. The cocycle
(\ref{cocyc}) is treated as an infinitesimal term of the
deformation
of the pairing $a\otimes b\mapsto \int_{{\cal O}}ab\,d\mu$
\cite{G2}.

In a similar way we consider an infinitesimal term of the
deformation of the pairing $A\otimes B\mapsto \Tr(A\circ B)$. For
this end, we use the relation
$$
\Tr_R\circ({\cal R}_{12} L_{1}{\cal R}_{21} L_{2} - L_{2}{\cal
R}_{12} L_{1}{\cal R}_{21}) = 0\,
$$
where matrix elements of the matrices $L_{1}$ and $L_{2}$ belong
to
$\End(V)$, the symbol $\circ$ stands for the product (\ref{l-tu})
in
this algebra and the operation $\Tr_R$ is applied to each matrix
element. The above relation holds true due to $\Tr_R\,
l_i^j=\delta_i^j$.

Then, expending the $R$-matrix and the $R$-trace into a series in
$\nu$:
$$
{\cal R} = I + \nu\,\textsf{r} + O(\nu^2)\,,\qquad
\Tr_R=\Tr+\nu\, b
+ O(\nu^2)\,,
$$
({\sf r} is given by (\ref{s-class})) we get the explicit form of
the operation $b\circ$ on the skew-symmetric subspace
$\wedge^2(\End(V))$
$$
b\circ(L_{1}\ot L_{2}-L_{2}\ot L_{1})=-\Tr\circ(
\textsf{r}_{12}L_{1}L_{2}+L_{1}\textsf{r}_{21}L_{2}-L_{2}
\textsf{r}_{12} L_{1}-L_{2}L_{1}\textsf{r}_{21})\,.
$$
Having thus defined the operation $b\circ$ on the basis elements,
we
directly get the general expression (see
(\ref{r-class}--(\ref{r-sup}) for notations) \be b\circ (A\otimes
B
- B\otimes A)=\Tr\circ (-r_-^{{\rm ad,ad}} (A\ot B)-r_+^{\rm
r,l}(A\ot B)+r_+^{\rm l,r}(A\ot B)),\quad A,B\in \End(V)\,.
\label{deftr} \ee

Thus, we have got the skew-symmetrized linear term of deformation
of
the pairing
$$
A\otimes B\mapsto \Tr(A\circ B)\,.
$$

\begin{proposition}
The quantity $\langle A,B\rangle = -b\circ (A\otimes B - B\otimes
A)$ is a cocycle on the Lie algebra $gl(m)$, i.e.
$$
\langle A,[B,C]\rangle+\langle B,[C,A]\rangle+\langle
C,[A,B]\rangle=0.
$$
It reduces to the Lie algebra $sl(m)$.
\end{proposition}

It is not difficult to write an explicit form of the cocycle
$\langle A,B\rangle$. Indeed, one can show that the second and
third
terms in the right hand side of (\ref{deftr}) give no
contribution
to this  cocycle and by using the cyclic property of the usual
trace we
get
$$
\langle A,B\rangle=\Tr\,([A,  B]\circ\sum_{\alpha >0} [X_\al,
X_{-\al}]) = \Tr\,([A,  B]\circ\sum_{\alpha >0} H_\alpha),
$$
where the sum is going over the set of all positive roots.

\section*{Appendix}
\setcounter{equation}{0}
\def\theequation{A.\arabic{equation}}
\medskip

In this Section we collect some facts and definitions on the
theory
of the $A_{k-1}$ series Hecke algebras $H_k(q)$, used in the main
text of the paper. For a detailed review of the subject the
reader
is referred to \cite{OP1}. Throughout this Section we use the
definitions and notations of that paper. At the end of the
Section
we give the proof of Proposition \ref{teo:3}, formulated in
Section
\ref{sec:2}.

By definition, a {Hecke algebra} of $A_{k-1}$ series is a unital
associative algebra $H_k(q)$ over a field $\Bbb K$ generated by
the
elements $\sigma_i$, $1\le i\le k-1$, subject to the following
commutation relations
$$
\begin{array}{lcl}
\sigma_i\sigma_{i+1}\sigma_i = \sigma_{i+1}
\sigma_i\sigma_{i+1} &\qquad&1\le i\le k-2\\
\rule{0pt}{5mm}
\sigma_i\sigma_j = \sigma_j\sigma_i & & |i-j| \ge 2\\
\rule{0pt}{5mm} \sigma_i^2 = 1_H - (q-q^{-1})\,\sigma_i & & 1\le
i\le k-1.
\end{array}
$$
Here $1_H$ is the unit of the algebra, $q\in \Bbb K$ is a nonzero
element of the ground field. Below we  assume $\Bbb K$ to be
the field of complex numbers $\Bbb C$ or the field of rational
functions ${\Bbb C}(q)$ of the formal variable $q$.

At a generic value of $q$ the Hecke algebra $H_p(q)$ is
semisimple
and isomorphic to the group algebra of the $k$-th order symmetric
group ${\Bbb K}[{\frak S}_k]$ \cite{We}. Therefore, being
considered
as the regular two-sided $H_k(q)$-module, the Hecke algebra
$H_k(q)$
can be presented as the direct sum of simple ideals (the
Wedderburn-Artin theorem)
$$
H_k(q) = \bigoplus_{\lambda\vdash k} M^{\lambda}
$$
labelled by partitions $\lambda$ of the integer $k$. Under the
left
(right) action of the Hecke algebra the submodules $M^{\lambda}$
are
reducible and can be further decomposed into the direct sum of
the
corresponding equivalent one-sided (left or right) submodules
$$
M^{\lambda} = \bigoplus_{a=1}^{d_\lambda}M^{(\lambda,a)}_{l(r)}
$$
where $d_\lambda$ is the number of the standard Young tableaux
$(\lambda,a)$ corresponding to the partition $\lambda$
\cite{Mac}.
The index $a$ enumerates standard tableaux in accordance with
some
ordering (say, lexicographical).

In each ideal $M^{\lambda}$ one can fix a linear basis of "matrix
units" $e^\lambda_{ab}$ with the multiplication table
$$
e^\lambda_{ab}\,e^\mu_{cd} =
\delta^{\lambda\mu}\delta_{bc}\,e^\mu_{ad}\,.
$$
A subset $e^\lambda_{ab}$, $1\le b\le d_\lambda$, (with a fixed
value of the first index) forms the basis of the right module
$M^{(\lambda,a)}_r$ while fixing the second index gives the basis
of
the left module $M^{(\lambda,b)}_l$.

The diagonal elements $e^\lambda_{aa}$ denoted shortly as
$e^\lambda_{a}$ form the set of primitive idempotents of the
Hecke
algebra $H_k(q)$. The idempotents $e^\lambda_a$ are explicitly
constructed as some polynomials in the Jucys-Murphy elements
${\cal
J}_p$, $1\le p\le k$, (see \cite{OP1} for details) which are
defined
by the iterative rule
$$
{\cal J}_1 = 1_H, \quad {\cal J}_{p+1} = \sigma_{p}\, {\cal
J}_{p}\,\sigma_{p}.
$$
The set of Jucys-Murphy elements form a basis of the maximal
commutative subalgebra of $H_k(q)$. An important property of
these
elements reads \be {\cal J}_p\, e^\lambda_a = e^\lambda_a\,{\cal
J}_p = j_p(\lambda,a)e^\lambda_a, \quad
j_p(\lambda,a)=q^{2(c_p-r_p)} \in \Bbb K. \label{j-calc} \ee Here
the positive integers $c_p$ and $r_p$ are the numbers of the
column
and the row of the Young tableau $(\lambda,a)$ which contain the
box
with integer $p$. Given below is a simple example for a Young
tableau of the partition $\lambda = (3,2,1)$
$$
\begin{tabular}{|c|c|c|}\hline
1&3&4\\ \hline 2& 6&\multicolumn{1}{c}{}\\ \cline{1-2}
5&\multicolumn{2}{c}{} \\ \cline{1-1}
\end{tabular}
\quad\Rightarrow\quad
\begin{array}{lcl}
j_1 = 1&\;&j_4 = q^4\\
j_2 = q^{-2}& & j_5 =q^{-4}\\
j_3 = q^2 & & j_6 = 1
\end{array}
$$

Any two idempotents $e^\lambda_a$ and $e^\lambda_b$ corresponding
to
the different tableaux $(\lambda,a)$ and $(\lambda,b)$ of a
partition $\lambda\vdash k$ can be transformed into each other by
the two sided action of some invertible elements of the Hecke
algebra $H_k(q)$ \cite{OP1} \be e^\lambda_a =
x^\lambda_{ab}\,e^\lambda_b\, y^\lambda_{ab}\,, \qquad
x^\lambda_{ab},\, y^\lambda_{ab}\in H_k(q). \label{id-conn} \ee

Consider now a Hecke symmetry $R:V^{\otimes 2}\to V^{\otimes 2}$,
and define a special representation $\rho_R$ of a Hecke algebra
$H_k(q)$ in the tensor product $V^{\otimes p}$, $p\ge k$, by the
rule \be
\begin{array}{l}
\rho_R(1_H) = {\rm id}_V^{\otimes p},\\
\rule{0pt}{5mm} \rho_R(\sigma_i) = {\rm id}_V^{\otimes
(i-1)}\otimes
R_{i}\otimes {\rm id}_V^{\otimes (p-i-1)}, \quad
1\le i\le k-1,\\
\rule{0pt}{5mm} \rho_R(xy) = \rho_R(x)\rho_R(y),\quad
\forall\,x,y\in H_k(q),
\end{array}
\label{r-rep-H} \ee recall, that $R_{i}:=R_{ii+1}$. The fact that
$\rho_R$ is a representation follows immediately from (\ref{YB})
and
(\ref{Hec}).

Let the bi-rank of $R$ be $(m|n)$, that is its HP series $P_-(t)$
is
of the form (\ref{p-}). Consider the partitions \be
\begin{array}{lcl}
\lambda_{m,n}:=((n+1)^{m+1})&\quad&\lambda_{m,n}\vdash
(m+1)(n+1)\\
\rule{0pt}{6mm} \lambda^-_{m,n}:= ((n+1)^{m},n)& &
\lambda^-_{m,n}\vdash mn+m+n\,.
\end{array}
\label{spec-d} \ee In graphic form the partition $\lambda_{m,n}$
is
represented by a rectangular diagram with $m+1$ rows of the
length
$n+1$, while the diagram of $\lambda^-_{m,n}$ is obtained from
the
former one by removing one box in the right lower corner of the
rectangle. Note, that $\lambda^-_{m,n}\in {\sf H}(m,n)$, while
the
partition $\lambda_{m,n}$ is the minimal one, not belonging to
the
hook ${\sf H}(m,n)$ (see Definition \ref{def:hook}).

Listed below are the properties of representations $\rho_R$ which
follow immediately from Proposition \ref{pro:dim-Vl}.
\begin{enumerate}
\item[i)] The images $E^\lambda_a = \rho_R(e^{\lambda}_a)
\not = 0$, $e^{\lambda}_a\in H_k(q)$, for all $2\le
k<(m+1)(n+1)$;
\item[ii)] the representation $\rho_R$ of $H_{(m+1)(n+1)}(q)$
possesses a kernel generated by
$$
\rho_R(e^{\lambda_{m,n}}_a)=0, \quad 1\le a \le
d_{\lambda_{m,n}}\,,
$$
while $\rho_R(e^\mu_a)\not=0$ for all $\mu\vdash (m+1)(n+1)$,
$\mu\not=\lambda_{m,n}$;
\item[iii)] for any integer $p\ge (m+1)(n+1)$ and for any
partition $\nu\vdash p$ one has
$$
\rho_R(e^\nu_a) = 0 \quad\Leftrightarrow \quad
\lambda_{m,n}\subset
\nu,
$$
where inclusion $\mu=(\mu_1,\mu_2,\dots)\subset \nu =
(\nu_1,\nu_2,\dots)$ means that $\mu_i\le \nu_i$ $\forall\, i$.
\end{enumerate}

\leftline{\bf Proof of Proposition \ref{teo:3}.} Let us denote
$p:=(m+1)(n+1)$ for more compact writing of the formulae below.
In
the Hecke algebra $H_p(q)$ we extract the Hecke subalgebra
$H_{p-1}(q)\subset H_p(q)$, generated by $\sigma_i\in H_p(q)$,
$1\le
i\le p-2$. Fix a standard Young tableau $(\lambda_{m,n},a)$ (see
definition (\ref{spec-d}) above) and consider the idempotents
$e^{\lambda^-_{m,n}}_{a^-} \in H_{p-1}(q)$ and
$e^{\lambda_{m,n}}_a\in H_p(q)$. Here the notation
$(\lambda^-_{m,n},a^-)$ refers to a special choice of the
corresponding Young tableau: it is properly included into the
Young
tableau $(\lambda_{m,n},a)$. In other words, the integers from
$1$
to $p-1$ occupy the same positions in the tableau
$(\lambda^-_{m,n},a^-)$ as they do in the tableau
$(\lambda_{m,n},a)$. Note, that since we consider the standard
Young
tableaux, the only possible position for the number $p$ is the
box
in the right lower corner of the rectangular tableaux
$(\lambda_{m,n},a)$.

Let us now apply the map $\rho_R:H_p(q)\rightarrow {\rm
End}(V^{\otimes n})$ to the relation (cf. \cite{OP1})
$$
e^{\lambda_{m,n}}_a = e^{\lambda^-_{m,n}}_{a^-}\, \frac{({\cal
J}_p
- q^{2(n+1)}1_H)}{(q^{2(n-m)} - q^{2(n+1)})}\, \frac{({\cal J}_p
-
q^{-2(m+1)}1_H)}{(q^{2(n-m)} - q^{-2(m+1)})}.
$$
Denoting $\rho_R({\cal J}_k):=J_k$ and taking into account the
item
ii) of the above properties of $\rho_R$ we get the identity
$$
0 = E^{\lambda^-_{m,n}}_{a^-}\, \frac{( J_p -
q^{2(n+1)}I)}{(q^{2(n-m)} - q^{2(n+1)})}\, \frac{(J_p -
q^{-2(m+1)}I)}{(q^{2(n-m)} - q^{-2(m+1)})}
$$
where $E^{\lambda^-_{m,n}}_{a^-}\not=0$ due to i), the letter $I$
stands for the identity operator on the space $V^{\otimes p}$.

We calculate the trace ${\sf tr}$ of the above identity in the
last ($p$-th) component of the tensor product $V^{\otimes p}$,
where
${\sf tr}$  coincides up to a factor with the categorical
$R$-trace
(\ref{trF})
$$
{\sf tr}(X):= Tr\,(C\cdot X).
$$
It is clear, that ${\sf tr}(I) = Tr\,C$ is the object we are
interested in.

Since the matrix $E^{\lambda_{m,n}^-}_{a^-}$ is a polynomial in
$J_k$ with $k<p$, it can be drawn out of the trace in the $p$-th
space and we come to \be 0 = E^{\lambda^-_{m,n}}_{a^-}\,{\sf
tr}_{(p)}\Bigl( J_p^2 - (q^{2(n+1)}+q^{-2(m+1)})J_p
+q^{2(n-m)}I\Bigr). \label{tr-id} \ee Consider separately the
traces
of the terms of the above identity. Introducing auxiliary
shorthand
notation $\omega:=q-q^{-1}$ we find
$$
{\sf tr}_{(p)}(J_p)  = {\sf
tr}_{(p)}\left(R_{p-1}J_{p-1}(R^{-1}_{p-1}+\omega I)\right)
 = \omega J_{p-1} + I_{p-1}\,{\sf tr}_{(p-1)}(J_{p-1}).
$$
In the above line of transformations we have used the iterative
definition of the Jucys-Murphy element, the Hecke condition for
$R$
and properties (\ref{trR}) and (\ref{trRC}) of ${\sf tr}$ listed
in
Section \ref{sec:2}. Since the trace in (\ref{tr-id}) is
multiplied
by the idempotent, we can replace the Jucys-Murhy element
$J_{p-1}$
by the corresponding "eigenvalue" $j_{p-1}$ defined in
(\ref{j-calc}) \be E^{\lambda^-_{m,n}}_{a^-}\,{\sf tr}_{(p)}(J_p)
=
E^{\lambda^-_{m,n}}_{a^-}\,\left(\omega j_{p-1} + {\sf
tr}_{(p-1)}
(J_{p-1})\right). \label{ind-1} \ee
To simplify the formulae, we
 omit the symbol of idempotent and perform all
calculations, bearing in mind the possibility to replace each
free
of trace Jucys-Murphy element $J_k$ by the corresponding number
$j_k$.

So, the calculation of ${\sf tr}_{(p)}(J_p)$ is completed by the
straightforward induction on the base of relation (\ref{ind-1})
$$
{\sf tr}_{(p)}(J_p) = \omega\sum_{k=1}^{p-1}j_k + {\sf tr}(I),
$$
where we have taken into account that $J_1 = I$ by definition.

Transform now the term with the second power of $J_p$:
\begin{eqnarray*}
{\sf tr}_{(p)}(J_p^2) &=& {\sf tr}_{(p)}
\left(R_{p-1}J_{p-1}R_{p-1}
(R^{-1}_{p-1}+\omega I)J_{p-1}R_{p-1} \right)  \\
&=&{\sf tr}_{(p)}\left(R_{p-1}J_{p-1}^2(R^{-1}_{p-1}+\omega
I)\right) + \omega \,J_{p-1}{\sf tr}_{(p)}
\left(J_p(R^{-1}_{p-1}+\omega I)\right)\\
&=&2\omega\, J_{p-1}^2+\omega^2\,J_{p-1}{\sf tr}_{(p)}(J_p) +
{\sf tr}_{(p-1)}(J_{p-1}^2)\\
&=&2\omega\, j_{p-1}^2+\omega^2\,j_{p-1}{\sf tr}_{(p)}(J_p) +
{\sf
tr}_{(p-1)}(J_{p-1}^2).
\end{eqnarray*}
Substituting the value of ${\sf tr}_{(p)}(J_p)$ we get the base
for
the inductive calculation
$$
{\sf tr}_{(p)}(J_p^2) = 2\omega\,j_{p-1}^2+\omega^2j_{p-1}\, {\sf
tr}(I) +\omega^3j_{p-1}\sum_{k=1}^{p-1} j_k+ {\sf
tr}_{(p-1)}(J_{p-1}^2).
$$
This immediately leads to the following expression
$$
{\sf tr}_{(p)}(J_p^2) = 2\omega \sum_{k=1}^{p-1}j_k^2+
\omega^3\sum_{k=1}^{p-1}j_k\sum_{s=1}^{k} j_s
+\Bigl(1+\omega^2\sum_{k=1}^{p-1}j_k \Bigr) \,{\sf tr}(I).
$$

Substituting all the calculated components into identity
(\ref{tr-id}) and taking into account that
$E^{\lambda^-_{m,n}}_{a^-}\not =0$, we find the following linear
equation for ${\sf tr}(I)$
$$
\alpha \,{\sf tr}(I) + \beta =0
$$
with
\begin{eqnarray*}
\alpha&=&1 + q^{2(n-m)} - q^{2(n+1)}-
q^{-2(m+1)} + \omega^2\sum_{k=1}^{p-1}j_k\\
\beta&=&\omega\left(2+\frac{\omega^2}{2}\right)
\sum_{k=1}^{p-1}j_k^2 + \frac{\omega^3}{2}
\Bigl(\sum_{k=1}^{p-1}j_k\Bigr)^2 - \omega\Bigl(q^{2(n+1)} +
q^{-2(m+1)}\Bigl) \sum_{k=1}^{p-1}j_k\,,
\end{eqnarray*}
where in finding the coefficient $\beta$ we used the identity
$$
\sum_{k=1}^{p-1}j_k\sum_{s=1}^{k}j_s =
\frac{1}{2}\sum_{k=1}^{p-1}j_k^2
+\frac{1}{2}\Bigl(\sum_{k=1}^{p-1}j_k \Bigr)^2.
$$
Taking into account the definition of $j_k$ and the form of the
diagram $\lambda_{m,n}^- = ((n+1)^m,n)$, we can easily calculate
the
sum of eigenvalues $j_k$
\begin{eqnarray*}
\sum_{k=1}^{p-1}j_k = \sum_{k=1}^{p-1}q^{2(c_k-r_k)} &=&
(1+q^2+\dots+q^{2n})(1+q^{-2}+\dots +q^{-2m}) - q^{2(n-m)}\\
&=& q^{n-m}(n+1)_q(m+1)_q - q^{2(n-m)}\,,
\end{eqnarray*}
and therefore
$$
\sum_{k=1}^{p-1}j_k^2 = \sum_{k=1}^{p-1}(q^2)^{2(c_k-r_k)} =
q^{2(n-m)}(n+1)_{q^2}(m+1)_{q^2} - q^{4(n-m)}\,.
$$
Now by a short calculation we simplify the coefficient $\alpha$
to
the form
$$
\alpha = -\omega^2\,q^{2(n-m)}.
$$
The transformation of $\beta$ is more involved though
straightforward too. In this way, one should use the following
identity
$$
k_{q^2}=\frac{q^{2k} - q^{-2k}}{q^2-q^{-2}} = \frac{(q^k -
q^{-k})}{(q-q^{-1})}\, \frac{(q^k + q^{-k})}{(q+q^{-1})} =
k_q\,\frac{q^k + q^{-k}}{2_q}.
$$
Omitting routine calculations we present the final result
$$
\beta = \omega^2\,q^{3(n-m)}(m-n)_q.
$$
So, we finally get
$$
{\sf tr}(I) = Tr\, C = -\frac{\beta}{\alpha} = q^{n-m}(m-n)_q.
$$
This completes the proof. \hfill \rule{6.5pt}{6.5pt}


\medskip

\begin{thebibliography}{IOP4}

\bibitem[AG]{AG} Akueson P., Gurevich D. {\em Dual
quasitriangular structures related to the Temperley-Lieb 
algebra} in: Lie groups and Lie algebras, Math. Appl., 
{\bf 433}, 1--16, Kluwer Acad. Publ., Dordrecht.

\bibitem[BR]{BR} Berele A., Regev A. {\em Hook Young diagrams
with applications to combinatorics and to representations of Lie
Superalgebras}, Adv. in Math., {\bf 64} (1987) 118--175.

\bibitem[BG]{BG} Braverman A., Gaistgory D. {\em The
Poincar\'e-Brkhoff-Witt theorem for quadratic algebras of Koszul
type}, J. Algebra {\bf 181} (1996) 315--328.

\bibitem[CP]{CP} Chari V., Pressley. A. {\em A guide to quantum
groups}, Cambridge University Press, Cambridge, 1994.

\bibitem[C]{C} Cherednik I. {\em Factorizing particles on a
half line, and root systems}, English translation in: Theor.
Math. Phys. {\bf 61} (1984) 977--983.

\bibitem[DGH]{DGH} Dabrowski L., Grosse H., Hajac P. {\em Strong
Connections and Chern-Connes pairing in the Hopf-Galois theory},
Comm. Math. Phys. {\bf 220} (2001) 301--331.

\bibitem[Da]{Da} Davydov A.A. {\it Totally positive sequences
and $R$-matrix quadratic algebras}, J. Math. Sci. {\bf 100}
(2001) 1871--1876.

\bibitem[DGG]{DGG} Delius G.W., Gardner C., Gould M.D. {\em The
structure of quantum Lie algebra for the classical series $B_l$,
$C_l$ and $D_l$}, J. Phys. A: Math. Gen. {\bf 31} (1998)
1995--2019.

\bibitem[D]{D} Donin J. {\em Double quantization on coadjoint
representations of semisimple Lie groups and their orbits},
{\sf math.QA/9909160}.

\bibitem[DGS]{DGS} Donin J., Gurevich D., Shnider S. {\em Double
quantization on some orbits in the coadjoint representation of
simple Lie groups}, Comm. Math. Phys. {\bf 204} (1999) 39--60.

\bibitem[DM]{DM} Donin J., Mudrov A. {\em Explicit equivariant
quantization on (co)adjoint orbits of $GL(n,\C)$}, Lett. Math.
Phys. {\bf 62} (2002) 17--32.

\bibitem[Dr]{Dr} Drinfeld V. {\em On quadratic
quasi-commutational relations in quasi-classical limit}, English
translation in: Selecta Math. Sovietica {\bf 11}, no. 4, (1992)
317--326.

\bibitem[DH]{DH} Dung N.P, Hai P.H. {\em On the Poincar\'e
series of quadratic algebras associated to Hecke symmetries}
Int. Math. Res. Not. {\bf 40} (2003) 2193--2203.

\bibitem[FP]{FP}
Faddeev L., Pyatov P. {\em The Differential Calculus on Quantum
Linear Groups}, Trans. Amer. Math. Soc., Ser. 2, {\bf 175} 
(1996) 35--47.

\bibitem[FH]{FH} Fulton W., Harris J. {\em Representation 
theory. A first course}, Springer-Verlag, NY, 1991.

\bibitem[FRT]{FRT} Reshetikhin N., Takhtadzhyan L., Faddeev L.
{\em Quantization of Lie groups and Lie algebras}, English
translation in: Leningrad Math. J. {\bf 1} (1990) 193--225.

\bibitem[GM]{GM} Gomez X., Majid S. {\em Braided Lie algebras and
bicovariant differential calculi over co-quasitriangular Hopf
algebras}, J. of Algebra, {bf 261} (2003) 334--388.

\bibitem[G1]{G1} Gurevich D. {\em Generalized translation
operators on Lie groups}, English translation in: Soviet J.
Contemporary Math. Anal. {\bf 18}, no. 4, (1983) 57--90.

\bibitem[G2]{G2} Gourevitch D. {\em Equation de Yang-Baxter et
quantification des cocycles}, C. R. Acad. Sci. Paris, {\bf 310}
(1990) 845--848.

\bibitem[G3]{G3} Gurevich D. {\em Algebraic aspects of the
Yang-Baxter equation}, English translation in: Leningrad Math. 
J. {\bf 2} (1991) 801 -- 828.

\bibitem[G4]{G4} Gurevich D. {\em Braided modules and reflection
equations}, Quantum Groups and Quantum spaces, Banach center
publications {\bf 40} (1995) 99--110.

\bibitem[GLS1]{GLS1}
Gurevich D., Leclercq R., Saponov P. {\em Traces in braided
categories}, J. Geom. Phys. {\bf 44} (2002) 251--278.

\bibitem[GLS2]{GLS2}
Gurevich D., Leclercq R., Saponov P., {\em q-Index on braided
noncommutative spheres}, J. Geom. Phys. {\bf 53} (2005) 
392--420.

\bibitem[GPS1]{GPS1}
Gurevich D., Pyatov P., Saponov P. {\em The Cayley-Hamilton
theorem for quantum matrix algebras of $GL(m|n)$ type}, English
translation in: St. Petersburg Math. J. {\bf 17} (2006) 
119--135.

\bibitem[GPS2]{GPS2}
Gurevich D., Pyatov P., Saponov P. {\em Quantum matrix algebras
of $GL(m|n)$ type: the structure of characteristic subalgebra
and its parametrization} (Russian), Teor. Mat. Fiz. {\bf 147}
(2006) 14--46.

\bibitem[GS1]{GS1}
Gurevich D., Saponov P. {\it Quantum line bundles via
Cayley-Hamilton identity}, J. Phys. A: Math. Gen. {\bf 34} (2001)
4553--4569.

\bibitem[GS2]{GS2}
Gurevich D., Saponov P. {\em On non-one-dimensional
representations of the reflection equation algebra}, English
translation in: Theor. Math. Phys. {\bf 139} (2004), 486--499.

\bibitem[GS3]{GS3}
Gurevich D., Saponov P. {\em Geometry of non-commutative orbits
related to Hecke symmetries}, Proceedings of International
Conference on Quantum Groups at Technion (J.Donin's memorial
volume will be published in Contemporary Mathematics).

\bibitem[GR]{GR} Gutt S., Rawnsley J. {\em Traces for star
products on symplectic manifols}, J. Geom. Phys. {\em 42} (2002)
12--18.

\bibitem[H]{H} Phung H.H. {\em Poincar\'e Series of Quantum
spaces Associated to Hecke Operators}, Acta Math. Vietnam 
{\bf 24} (1999) 235--246.

\bibitem[IOP]{IOP} Isaev A.,  Ogievetsky O., Pyatov P.
{\em On quantum matrix algebras satisfying the
Cayley-Hamilton-Newton identities}, J. Phys. A: Math. Gen. 
{\bf 32} (1999) L115--L121.

\bibitem[IP]{IP}
Isaev A.,  Pyatov P. {\em Covariant differential complexes on
quantum linear groups}, J. Phys. A: Math. Gen. {\bf 28} (1995)
2227-2246.

\bibitem[KT]{KT} Khoroshkin S., Tolstoy V. {\em Universal
$R$-matrix for quantized (super-)algebras}, Comm. Math. Phys. 
{\bf 141} (1991) 599--617.

\bibitem[K]{K} Kulish P.P. {\em Representations of
$q$-Minkowski space algebra}, St. Petersburg Math. J. {\bf 6}
(1995) 365--374.

\bibitem[KS]{KS} Kulish P., Sklyanin E. {\em Algebraic structure
related to the reflection equation}, J. Phys. A: Math. Gen 
{\bf 25} (1992) 5963--5975.

\bibitem[LW]{LW} Lu J.-H.,Weinstein A. {\em Poisson Lie groups,
dressing transformations, and Bruhat decompositions}, J.
Differential Geom. {\bf 31} (1990) 501--526.

\bibitem[LS]{LS} Lyubashenko V., Sudbery A. {\em Quantum
supergroups of $GL(m|n)$ type: differential forms, Koszul
complexes, and Berezinians}, Duke Math. J. {\bf 90} (1997) 
1--62.

\bibitem[Mac]{Mac} Macdonald I. G. {\em Symmetric Functions and
Hall Polynomials (Oxford Mathematical Monographs)}, Oxford
Science Publications, 1995.

\bibitem[M]{M} Majid S. {\em Foundations of quantum group
theory}, Cambridge University Press, Cambridge, 1995.

\bibitem[Mu1]{Mu1} Mudrov A. {\em Characters of
$U_q(gl(m))$-reflection equation algebra}, Lett. Math. Phys 
{\bf 60} (2002) 283--291.

\bibitem[Mu2]{Mu2} Mudrov A. {\em On quantization of
Semenov-Tian-Shansky Poisson bracket on simple algebraic groups}
{\sf math.QA/0412360}.

\bibitem[O]{O} Ogievetsky O. {\em Uses of Quantum Spaces},
Contemp. Math. {\bf 294} (2002) 161--232, AMS, Providence, RI.

\bibitem[OP1]{OP1} Ogievetsky O. and Pyatov P. {\em Lecture on
Hecke algebras}, in Proc. of the International School 
"Symmetries and Integrable Systems" Dubna, Russia, June 8-11,
1999. JINR, Dubna, D2, 5-2000-218, 39--88; Preprints
CPT-2000/P.4076 and MPI 01-40.

\bibitem[OP2]{OP2} Ogievetsky O., Pyatov P. {\em
Orthogonal and Symplectic Quantum Matrix Algebras and
Cayley-Hamilton Theorem for them}, {\sf math.QA/0511618}.

\bibitem[P]{P} Podles P. {\em Quantum spheres}, 
Lett. Math. Phys. {\bf 14} (1987) 193--202.

\bibitem[PP]{PP} Polishchuk A., Positselski L. {\em Quadratic
Algebras}, University Lecture Series, {\bf 37}, AMS, Providence,
RI.

\bibitem[R]{R} Rieffel M. {\em Deformation quantization of
Heisenberg manifolds}, Comm. Math. Phys. {\bf 122} (1989)
531--562.

\bibitem[S]{S} Saponov P. {\em The Weyl approach to the
representation theory of reflection equation algebra}, J. Phys.
A: Math. Gen. {\bf 37} (2004) 5021--5046.

\bibitem[STS]{STS} Semenov-Tian-Shansky {\em Dressing
transformations and Poisson group actions}, Publ. Res. Ins. Math.
Sci., {\bf 21} (1985) 1237-1260.

\bibitem[Sh]{Sh} Sheu A. {\em Quantization of the Poisson $SU(2)$
and Its Poisson Homogeneous Space-The 2-Sphere}, Comm. Math.
Phys. {\bf 135} (1991) 217--232.

\bibitem[St]{St} Stembridge J.R. {\em A Characterization of
Supersymmetric Polynomials}, J. Algebra {\bf 95} (1985) 
439--444.

\bibitem[T]{T} Turaev V. {\em Quantum invariant of knots and
3-manifolds} W. de Gruyter, Berlin, 1994.

\bibitem[We]{We} Wenzl H. {\em  Hecke algebras of type $A_n$ and
subfactors}, Invent. Math. {\bf 92} (1988) 349--383.

\bibitem[Wo]{Wo} Woronowicz S. {\em Differential Calculus on
Compact Matrix Pseudogroups (Quantum groups)}, Comm. Math. Phys.
{\bf 122} (1989) 125--170.

\bibitem[Z]{Z} Zhang R.B. {\em Structure and representations of
the quantum general Linear supergroup}, Comm. Math. Phys. 
{\bf 195} (1998) 525--547.


\end{thebibliography}
\end{document}